\documentclass[12pt,psfig]{article}
\usepackage{amsthm,amssymb,amsmath,amsfonts,amscd}
\usepackage{graphicx,epsfig,latexsym}
\usepackage{cite,calc,xcolor,subfigmat,mathrsfs,dsfont}

\newcommand{\F}{\mathscr{F}}
\newcommand{\LST}{\mathscr{L}}

\newcommand{\Y}{\mathcal{Y}}

\newcommand{\FS}{\mathscr{F}}
\newcommand{\Sl}{\mathfrak{sl}_2}
\newcommand{\T}{\mathfrak{T}}

\newcommand{\f}{\mathbb{F}}

\newcommand{\N}{\mathbb{N}}

\newcommand{\Span}{{\rm span}}

\newcommand{\ad}{{\rm ad}}

\newcommand{\rank}{{\rm rank}}
\newcommand{\0}{\mathbf{0}}

\newcommand{\C}{\mathcal{C}}
\newcommand{\m}{\mathbf{n}}

\newcommand{\X}{\mathbb{F}}
\newcommand{\XS}{\mathcal{X}}
\newcommand{\IM}{{\rm Im}\,}
\newcommand{\ie}{{\em i.e.,} }
\newcommand{\eg}{{\em e.g.,} }

\newcommand{\y}{\Delta}

\newcommand{\ddx}{\frac{\partial}{\partial x}}
\newcommand{\ddy}{\frac{\partial}{\partial y}}
\newcommand{\ddz}{\frac{\partial}{\partial z}}

\newtheorem{thm}{Theorem}[section]

\newtheorem{lem}[thm]{Lemma}
\newtheorem{prop}[thm]{Proposition}
\theoremstyle{definition}

\newtheorem{exm}[thm]{Example}
\newtheorem{rem}[thm]{Remark}
\newtheorem{notation}[thm]{Notation}
\numberwithin{equation}{section}
\def\be {\begin{equation}}
\def\ee {\end{equation}}
\def\ba {\begin{eqnarray}}
\def\ea {\end{eqnarray}}
\def\bes {\begin{equation*}}
\def\ees {\end{equation*}}
\def\bas {\begin{eqnarray*}}
\def\eas {\end{eqnarray*}}
\def\bpr {\begin{proof}}
\def\epr {\end{proof}}

\begin{document}
\baselineskip=18pt
\renewcommand {\thefootnote}{\dag}
\renewcommand {\thefootnote}{\ddag}
\renewcommand {\thefootnote}{ }

\pagestyle{empty}

\begin{center}
\leftline{}
\vspace{-0.500 in}
{\Large \bf Volume-preserving normal forms of Hopf-zero singularity} \\ [0.3in]

{\large Majid Gazor$^{*},$ Fahimeh Mokhtari}
\footnote{$^*\,$Corresponding author. Phone: (98-311) 3913634; Fax: (98-311) 3912602; Email: mgazor@cc.iut.ac.ir}

\vspace{0.105in} {\small {\em Department of Mathematical Sciences,
Isfahan University of Technology
\\[-0.5ex]
Isfahan 84156-83111, Iran }}

\today

\vspace{0.05in}

\noindent

\end{center}

\vspace{-0.10in}
\baselineskip=15pt
\noindent \rule{6.5in}{0.012in}

\vspace{-0.1in} \noindent

\begin{abstract}

A practical method is described for computing the unique generator of the algebra of first integrals associated with a large class of Hopf-zero singularity. The set of all volume-preserving classical normal forms of this singularity is introduced via a Lie algebra description.
This is a maximal vector space of classical normal forms with first integral; this is whence our approach works.
Systems with a non-zero condition on their quadratic parts are considered. The algebra of all first integrals for any such system has a unique
(modulo scalar multiplication) generator. The infinite level volume-preserving parametric normal forms of any non-degenerate
perturbation within the Lie algebra of any such system is computed, where it can have rich dynamics.
The associated unique generator of the algebra of first integrals are derived. The symmetry group of the infinite level normal forms are also discussed.
Some necessary formulas are derived and applied to appropriately modified R\"{o}ssler  and generalized Kuramoto--Sivashinsky equations to demonstrate the applicability of our theoretical results. An approach (introduced by Iooss and Lombardi) is applied to find an optimal truncation for the first level normal forms of these examples with exponentially small remainders. The numerically suggested radius of convergence (for the first integral) associated with a hypernormalization step is discussed for the truncated first level normal forms of the examples. This is achieved by an efficient implementation of the results using Maple.

\vspace{0.10in} \noindent {\it Keywords}: \ Volume-preserving; First integral; Normal form; Hopf-zero singularity; Lie algebra.

\vspace{0.10in} \noindent {\it 2010 Mathematics Subject Classification}:\, 34C20, 34A34, 16W50, 68U99.
\end{abstract}
\vspace{-0.2in}
\noindent \rule{6.5in}{0.012in}

\vspace{0.2in}

\section{Introduction}

A Lie algebraic structure for a family of vector fields is important in the sense that the family is invariant under a group of permissible
transformations generated by the Lie algebra. Thus, certain dynamical properties are preserved within the family. Therefore, in the study
of dynamical systems, it is fundamentally useful to recognize nonlinear families that admit a Lie algebra structure. Once this is accomplished, the transformation group makes a classification for the vector fields. This classification can be computed through the \emph{infinite level} (\emph{simplest}, \emph{unique}) normal form theory. This gives rise to the study of singular differential systems which has significantly contributed to the bifurcation and
stability analysis of many systems with complex dynamics. A good practical knowledge on the Lie algebraic structure of any
singularity is necessary for an efficient infinite level normal form computation. There are many results on the simplest normal forms of planar singularities. However, there are considerably less research results on the simplest normal forms of the three-dimensional state space
singularities, where the system is already on its center manifold; see \cite{ChenHopfZ03,ChenHopfZ,Algaba3Z,YuHopfZero,YuYaun3Zero}. The reason is that the Lie
algebraic structure generated by three dimensional state space is more complicated than that of the planar state space.
Findings on the Lie subalgebras of a singularity contribute to the understanding of the dynamical properties invariant
under the transformation groups associated with the Lie subalgebras. We introduce a maximal Lie algebra of Hopf-zero normal forms with a first integral which are derived via a \(\Sl\)-representation in \cite{GazorMokhtari}. By further investigation, it turns out that this family is the set of all volume-preserving Hopf-zero normal forms.

\pagestyle{myheadings}
\markright{{\footnotesize {\it  M. Gazor and F. Mokhtari
\hspace{2.65in} {\it Volume-preserving normal forms }}}}

A practical method for finding the generators of the algebra generated by first integrals of a given differential system is
an indispensable subject. Our approach displays such a fruitful method that works well for a large class of Hopf-zero singular systems.
The idea is to transform the system into the classical normal form and check if it is a volume-preserving normal form system. Then, our formulas derive the unique generator for the algebra of its first integrals. Therefore, this method is applicable for all Hopf-zero singular systems that their classical normal form is divergent free.

Our proposal for normalization of a Hopf-zero singularity is as follows. A vector field is first transformed into the classical normal form and then, we consider its conservative--nonconservative decomposition; see \cite{BaidSand91,baidersanders,PalacSIAM,PalacNonl,PalacChaos} for some relevant results. Based on this observation, an appropriate hypernormalization approach is applied. This paper deals with hypernormalization for the cases when the nonconservative part is zero. For the zero conservative part, the hypernormalization follows from the method introduced in \cite{GazorMokhtariEul}. Finally, for the cases of both nonzero conservative and nonconservative parts, the normal form computation in \cite{GazorMokhtari} must be implemented. The hypernormalization steps associated with either of these are independent and Maple programs are developed for their implementations.

For preserving the divergent-free property of the normal form systems, the transformation group (or instead, transformation group generators) must be defined appropriately. The normal form study of three-dimensional volume-preserving vector fields have also been considered by Mezi\'{c} and Wiggins \cite{Wiggins3DMez}. Their results are consistent with our claims and this suggests that divergent free Hopf-zero singularity may appear in models of incompressible fluid flows; see Remark \ref{RemWigginsMezic}. For some relevant results associated with volume preserving maps see \cite{Bazzani}.
Here, we appreciate the editor-in-chief careful readying of our manuscript and fruitful comments. The present paper focuses on three-dimensional divergence-free vector fields. Such vector fields arise in incompressible hydrodynamics (Lagrangian motion of particles) and in studies of magnetic field lines in three dimensions; \eg see \cite{Wiggins3DMez,Wiggins,LauFinn}. Near stagnation points of a flow or magnetic nulls the vector field is degenerate and normal form theory provides valuable information about the behavior of the system.

In this paper we are concerned with computing the simplest normal form of the divergent free systems given by
\ba\label{Eq1}
\left\{
  \begin{array}{ll}
\dot{x}= \sum_{} (k-l+1)a^l_k x^{l+1} {(y^2+z^2)}^{k-l}, &
\\
\dot{y}= -\sum_{} \frac{(l+1)}{2}a^l_kx^{l}y {(y^2+z^2)}^{k-l}+\sum_{} b^m_nx^{m}z(y^2+z^2)^{n-m},&
\\
\dot{z}= -\sum_{} \frac{(l+1)}{2}a^l_kx^{l}z {(y^2+z^2)}^{k-l}- \sum_{} b^m_nx^{m}y {(y^2+z^2)}^{n-m}, &
\end{array}
\right.
\ea where
\bes -1\leq l\leq k, 0\leq k, \quad 0\leq m\leq n,\quad (x,y,z)\in \mathbb{R}^3, \quad a^l_k, b^m_n\in \mathbb{R}, \quad a^{-1}_0\neq 0, a^0_0=0,\ees
and \( b^0_0=1.\) This system is the classical normal form of a large class of Hopf-zero singularity systems with a first integral. For instance Guckenheimer and Holmes \cite[Equations 7.4.5 for \(a:=1\)]{Guckenheimer} studied a generic example from this family. The formal function
\ba\label{Eq2}
f(x, y, z)&:=& \sum a^l_k x^{l+1} {(y^2+z^2)}^{k-l+1} \qquad \hbox{ (where } -1\leq l\leq k\hbox{ and } l+k\geq 1)
\ea is a first integral for the system \eqref{Eq1}. The idea is to analyze the dynamics from the infinite level normal form
and its first integral rather than working with the more complicated system \eqref{Eq1} and first integral \eqref{Eq2}.

The rest of this paper is organized as follows. Section \ref{secDec} discusses the notations and algebraic structures associated with the volume-preserving
classical Hopf-zero normal form with a first integral. We obtain the infinite level volume-preserving normal forms and associated first integrals
in Section \ref{SNF}. Furthermore, a new technique is demonstrated in which it enables one to change certain coefficients into non-algebraic numbers in order to facilitate elimination of certain terms. Then after some normal form computation, one changes them back into the desired numbers; see Remark \ref{arbitA}.
In Section \ref{PNF}, under some technical conditions, we prove that there exist invertible transformations
sending any non-degenerate perturbation of the system \eqref{Eq1} into a truncated infinite level parametric normal form (in cylindrical coordinates)
\ba\label{InftIntr}
\left\{
  \begin{array}{ll}
\dot{x}&= \rho^2\pm  x^{p+1}+\sum^r_{i=1} x^{n_i+1}\mu_i+\sum^{N}_{k=p+1} \alpha_{k} x^{k+1},\\
\dot{\rho}&=\mp\frac{(p+1)}{2} x^{p}\rho-\sum^r_{i=1} \frac{n_i+1}{2}x^{n_i}\rho\mu_i-\rho\sum^{N}_{k=p+1}\frac{(k+1)\alpha_{k}}{2}x^{k}, \\
\dot{\theta}&=1+ \sum^{r+s}_{i=r+1} x^{m_{i-r}}\mu_i+ \sum^{N}_{l=q} \beta_{l}x^l,
  \end{array}
\right.
\ea for some \bes r, s, p, q\geq 1, N-1\geq n_i\geq -1, \hbox{ and } N-1\geq m_i\geq 1;\ees see Equations \eqref{nj} and \eqref{mj}.
Here, \(\alpha_{k}=0\) for \(k\equiv p-1 \mod{2(p+1)},\) while
\bes \beta_{l}=0 \hbox{ for } l\equiv -1\mod{2(p+1)} \hbox{ and for } l\equiv p+q\mod{2(p+1)}. \ees
Furthermore, the function
\bes
f(x,\rho)=\rho^2\left( \frac{1}{2} \rho^2\pm x^{p+1}+\sum^r_{i=1} x^{n_i+1}\mu_i+\sum^{N}_{k=p}\alpha_{k} x^{k+1}\right)
\ees is the unique generator (modulo scalar multiplication) for the algebra of first integrals for system \eqref{InftIntr}; see Theorem \ref{InftPNF}. In Section \ref{Examples} we provide the necessary formulas for a Hopf-zero singular system whose cubic-truncated first level normal form is governed by Equation \eqref{Eq1}. Using these formulas, we first apply our approach to a modified R\"{o}ssler equation in Example \ref{RosselerExample} and then to a symmetric system chosen from a commonly used generalized Kuramoto--Sivashinsky equation in Example \ref{KSExample}. Section \ref{IoosSection} is devoted to find an optimal truncation for the first level normal forms. We restate some results from Iooss and Lombardi \cite{IoosLombardi} in an specific form we need. We apply them to the examples of Section \ref{Examples} which are presented in Examples \ref{RosselerExample6} and \ref{KSExample6}. All formulas are implemented using Maple. Finally, the paper is concluded in Section \ref{SecConv} with a convergence analysis based on numerically computed first integrals for normal forms associated with Examples \ref{RosselerExample} and \ref{KSExample}. These are respectively presented in Examples \ref{RosselerExample7} and \ref{KSExample7}.

\section{Volume-preserving conservative vector fields } \label{secDec}

An original motivation was to consider normal forms of conservative Hopf-zero singularity systems. Conservative systems have many applications in real life problems and Hamiltonian systems are among the most prominent examples.
There are substantial contributions on Hamiltonian systems in the literature. However, the Hamiltonian
structure requires an even dimensionality of the state space. Thus, an eminent alternative for conservative Hopf-zero singular systems on center manifold,
\ie a three dimensional state space, is to consider vector fields with a first integral; see \cite{LlibreWalcher} for some relevant results. A normal form computation may not destroy certain symmetric structures (\eg volume-preserving, Hamiltonian, etc.) when the transformation group is generated by an appropriate Lie algebra that preserve symmetric (volume-preserving, Hamiltonian) vector fields.
However, the set of all classical normal Hopf-zero singularities with a first integral neither is closed under the Lie bracket nor is a vector space.
Therefore, we do not consider the {\em set} of all vector fields with a first integral. Instead, we introduce a maximal vector
space (that is also a maximal Lie algebra) of such vector fields. The systems introduced in this paper are derived from a \(\Sl\)-representation for the classical normal forms of Hopf-zero singularity; see \cite{GazorMokhtari}. Here, the basic ideas stem from the lessons that the first author learned from Professor Jan A. Sanders in his 2010 summer visit of Vrije university. Here, we discuss the algebraic structures and dynamics properties in details. The presented algebraic structures are necessary for the normal form computation in the following sections. For terminologies and background used in this paper see \cite{GazorYuSpec,MurdBook,Murd04,Mokhtari}.

We denote
\ba
F^l_k&:= &x^{l} {(y^2+z^2)}^{k-l}\Big((k-l+1) x\ddx-
\frac{(l+1)}{2} y\ddy-\frac{(l+1)}{2} z\ddz\Big),\quad   (-1\leq{l}\leq{k})\;\;\;
\\ \nonumber
\Theta^l_k&:=& x^{l}(y^2+z^2)^{k-l}z\ddy- x^{l} {(y^2+z^2)}^{k-l}y\ddz,\qquad
\qquad\qquad\qquad\qquad\,\;\;
(0\leq{l}\leq{k}).\quad\;\;\;
\ea
The vector fields \(F^l_k\) and \(\Theta^l_k\) can also be represented in cylindrical coordinates, but
many first integrals for \(F^l_k\) only appear in terms of \(x, y, z\); see Lemma \ref{Lem1int}.
For a \(0\neq a^{-1}_0\in \mathbb{R},\) let
\be\label{Lie}
\LST:= \Span\left\{\Theta^0_0+a^{-1}_0 F^{-1}_0+ \sum a^l_k F^l_k+\sum b^m_n
\Theta^m_n\,|\, k, n\geq 1, a^l_k, b^m_n\in \mathbb{R}\right\}, \ee where  \(-1\leq l\leq k\)  and \(0\leq m\leq n.\)
These notations provide a tool to use a similar approach to the method developed
by Baider and Sanders \cite{BaidSand91,baidersanders}; see also \cite{GazorMoazeni}.
The space \(\LST\) is a Lie algebra by the Lie bracket \([v,w]= vw-wv\) for any \(v,w\in \LST.\)
Let
\bes\T:= \Span\left\{\Theta^0_0+\sum b^l_k \Theta^l_k\,|\, 0\leq l\leq k, 1\leq k, b^l_k\in \mathbb{R}\right\}.\ees
This is
equivalent to the space of phase components in cylindrical coordinates.
The set of all formal first integrals for \(v\in \LST\) is a subalgebra of formal power
series in terms of \(x, y, z.\)
Since the formal power series is a Noetherian ring, it is finitely generated.
We denote by \(\langle \,p_1, p_2, \ldots, p_k\rangle\) the algebra generated by
\(p_1, p_2, \ldots, p_k\in \mathbb{R}[[x, y, z]],\) unless otherwise is stated.
\begin{lem} \label{Lem1int} For any \(l\) and \(k,\)
let \(\FS^l_k\) be the algebra of formal
first integrals for \(F^l_k.\) Then we have:
\begin{itemize}
  \item\label{Part1} The algebra of first integrals for any \(0\neq v\in \T\) is
\bas
\left\langle x, y^2+z^2\right\rangle.
\eas
  \item The algebra \(\FS^l_k\) is generated by monomials that their \(x\)-degree is \((l+1)\) and their \((y,z)\)-degree is \(2(k-l+1),\) \ie
\bas
\FS^l_k=\left\langle x^{l+1}y^{2i}z^{2(k-l+1)-2i}\,\big|\, i=0, 1, \ldots, k-l+1 \right\rangle.
\eas
  \item For any \(0\neq a\in \mathbb{R},\)
\bas
\left\langle x^{l+1}(y^2+z^2)^{k-l+1}\right\rangle
\eas is the algebra of first integrals for \(\Theta^0_0+ aF^l_k.\)
\end{itemize}
\end{lem}

\bpr Let \(g\) be a formal first integral for \(\Theta^m_n.\) Then, \(z\ddy g-y\ddz g=0\) and
\bes g=g\left(x, y^2+z^2\right)\in \left\langle x, y^2+z^2\right\rangle.\ees The second part is a straightforward computation.
Since every formal first integral for \(\Theta^0_0+ aF^l_k\) must be a first integral for \(\Theta^0_0\) (\eg see \cite[Proposition 3]{LlibreWalcher}),
any first integral for \(\Theta^0_0+ aF^l_k\) belongs to
\(\langle x, y^2+z^2\rangle\cap \FS^l_k.\) This completes the proof.
\epr
The above lemma suggests the following proposition.
\begin{prop}\label{1integ}
Let \(v= \Theta^0_0+a^{-1}_0 F^{-1}_0+ \sum a^l_k F^l_k+\sum b^l_k \Theta^l_k\in \LST,\) where \(a^{-1}_0\neq 0.\)
Then, there exists a unique formal first integral \(f\) (modulo scalar multiplications) such that the algebra of first integrals for
\(v\) is \(\langle f\rangle.\)
\end{prop}

\bpr Define \bes f:= \int_0^1 \left\langle v\left( tx, t^{\frac{1}{2}}y, t^{\frac{1}{2}}z \right), \left(y^2+z^2, -2xy, -2xz\right)\right\rangle\, d t,\ees
where \(\langle\cdot,\cdot\rangle\) denotes the usual
inner product on \(\mathbb{R}^3.\) Then,
\be f=a^{-1}_0(y^2+z^2)^{2}+\sum a^l_kx^{l+1}(y^2+z^2)^{k-l+1}.\ee
The formal function \(f\) is a formal first integral for \(v\) because \(v(f)=0.\)
By \cite[Proposition 3]{LlibreWalcher} any first integral for \(v\) lies in \(\langle x, y^2+z^2\rangle\) and is also a first integral for
\(w:=a^{-1}_0 F^{-1}_0+ \sum a^l_k F^l_k.\) Let \(g(x, y^2+z^2)\) be such that
\bes w(g)=0 \hbox{ and } g= g_N+ \cdots,\ees
where \(g_N\) denotes the nonzero homogenous polynomial component of \(g\) with the least degree. Since \(w(g)=0,\) the least degree of the expansion for
\(w(g)\) must be zero, \ie \(F^{-1}_0(g_N)=0.\)
By Lemma \ref{Lem1int}, \(g_N\in \langle (y^2+z^2)^2\rangle\) and there exist a natural number \(k\) and a real number \(a_1\in \mathbb{R}\) such that
\bes N= 4k_1\hbox{ and } g_N= a_1\left(y^2+z^2\right)^{2k_1}.\ees
Now let \(\tilde{g}_1:= g- f^{k_1}.\) Then, \(\tilde{g}_1\) is also a first integral for \(w\) whose nonzero homogenous polynomials (monomials) in its power series expansion have degrees strictly greater than \(N.\) An induction argument proves that there exist real numbers \(a_i\in \mathbb{R}\) and natural numbers \(k_i\) (for \(i=2, 3, \ldots\)) such that \bes\tilde{g}_i:= g- \sum^{i-1}_{j=1} a_jf^{k_j}.\ees The proof is complete since
\(\tilde{g}_i\) is convergent to zero with respect to the filtration topology.
\epr

\begin{rem}
Any \(v\in \LST\) has a first integral and is volume-preserving, \ie \({\rm div}(v)=0. \) However, \(\LST\) is not
the {\it set} of all vector fields with first integral.
Indeed, any function of the form \(g(y^2+z^2)\) is a first integral for
\bes w_l:=F^l_l+\frac{(l+1)}{2}x^l\Big(2x\ddx+ y\ddy+ z\ddz\Big)= (l+2)x^{l+1}\ddx,\ees for any \(l\in \N.\) However, the vector fields
\(w_l-F^l_l\)
and
\bas
[w_5, F^2_3]&=& -21 x^7\left(y^2+z^2\right)\left(2x\ddx+ y\ddy+z\ddz\right)
\eas do not have any first integral and are not volume-preserving; see \cite[Theorem 2.3]{GazorMokhtariEul}.
This implies that the set of all vector fields with a first integral is not closed under the Lie bracket and it is not a vector space.
Indeed, any classical normal form with Hopf-zero singularity is uniquely decomposed into a divergent free vector field (with a first integral)
from \(\LST\) and a vector field with non-zero divergent (without any first integral); see \cite{GazorMokhtari}. This implies that \(\LST\) is the set of all volume-preserving Hopf-zero classical normal forms.
\end{rem}

\begin{rem}\label{RemWigginsMezic}
Mezi\'{c} and Wiggins \cite{Wiggins3DMez} considered a class of three-dimensional systems associated with incompressible (volume-preserving) fluid flows. They considered such vector fields when they admit a one-parameter spatial volume-preserving symmetry group. They proved that there exists a local change of variables such that it sends the system into an analytic normal form. The significance of their normal form is that the evolution of two variables is governed by a one-degree of freedom Hamiltonian system; while the evolution of the third variable depends only on the first two variables. When the original system is autonomous, the three-dimensional normal form system has a first integral; see \cite[Theorem 2.1]{Wiggins3DMez} and \cite[Theorem 2.66]{Olver}. Our normal forms are consistent with their results. Notice that the system \eqref{Eq1} admits a one-parameter spatial volume-preserving symmetry group whose infinitesimal generator is given by \(z\ddy-y\ddz\). Here both the divergent free property and the rotational (a one-parameter) symmetry group are discovered after the first level normal forms are computed. However, the considered family in \cite{Wiggins3DMez} already must have these symmetries in order to be transformed into normal forms. This suggests that analytic systems governed by Equation \eqref{Eq1} may be derived from practical models in incompressible fluids.
\end{rem}

The following lemma portrays the structure constants involved in this paper.

\begin{lem} The following equations hold true.
\bas\label{aa}
\left[F^l_k, F^{m}_{n}\right]&=&
\big((m+1)(k+2)-(l+1)(n+2)\big)F^{l+m}_{k+n}, \;\hbox{ for } \;\; -1\leq l\leq k, -1\leq m\leq n,
\\
\left[F_k^l,\Theta_n^m\right]&=&\big(m(k+2)-n(l+1)\big)\Theta^{l+m}_{ k+n},\qquad\qquad\quad\, \hbox{ for }\;\;\,-1\leq l\leq k, 0\leq m\leq n,
\\
\left[\Theta_k^l,\Theta_n^m\right]&=&0, \qquad\qquad\qquad\qquad\qquad\qquad\qquad\qquad\; \hbox{ for } \quad\quad 0\leq l\leq k, 0\leq m\leq n.
\eas
\end{lem}
The space \(\T\) is a nontrivial Lie ideal (and a trivial Lie subalgebra) for \(\LST.\) The above lemma implies that the quotient
Lie algebra $\LST/\T$ is Lie-isomorphic to a proper Lie subalgebra of a one-degree of freedom Hamiltonian vector fields;
see \cite[Theorem 3.7]{baidersanders} and \cite[Theorems 2.5 and 2.7]{GazorMokhtari} for more detailed discussions.

\begin{notation}
Throughout this paper, we use Pochhammer \(k\)-symbol notation, that is,
\bes
(a)^k_{b}:= a (a+b)(a+2b)\ldots \big(a+(k-1)b\big),
\ees for any natural number \(k\) and real number \(b.\)
\end{notation}

\section{The infinite level normal forms} \label{SNF}

In this section, we obtain the simplest normal form of Hopf-zero systems given by Equation \eqref{Eq1}.
Let \(v= \sum^\infty_{k=0}v_k\) be a Lie-graded expansion of \(v\in \LST.\)
Define \(d^{k,1}: \LST_k\rightarrow \LST_k\) by \(d^{k,1}(Y_k)= [Y_k, v_0]\) and then inductively define the maps
\(d^{k, n}: \LST_k\times \ker d^{k,n-1}\rightarrow \LST_k\) by
\bes
d^{k,n}\left(Y_{k}, Y_{k-1}, \ldots , Y_{k-n+1}\right):= \sum^{n-1}_{i=0} \left[Y_{k-i}, v_i\right], \quad \hbox{ for any } k\geq n.\ees For \(k<n,\) define
\(d^{k,n}:= d^{n,n}.\) A normal form style is a rule on how to choose a unique complement space \(\C^{k,n}\) so that \(\IM d^{k,n}\oplus \C^{k,n}=\LST_k.\)
Then, for a given vector field \(v,\) a graded Lie algebra structure and a normal form style, there exists a formal invertible transformation
that transforms \(v\) into its \(n\)-th (infinite or simplest) level normal form \(w\) where \(w= \sum^\infty_{k=0} w_k\) and \(w_k\in \C^{k,n}\)
(\(w_k\in \C^{k,k}\)) for all \(k\); see \cite{baiderchurch,GazorYuSpec} for details.

In the literature \(\Sl\)-style normal form has been applied only as a {\it first level style} (see \cite{Murd08,Murd04,MurdBook} for more details on {\it normal form styles}), while we use it as a {\it second level style} in the following lemma. In other words, the \(\Sl\)-representation is here applied to the first level normal forms (of Hopf-zero singularity) and then, \(\Sl\)-style is applied to the second level normal form computation. Derivation of this family through \(\Sl\)-representation is beyond the scope of this paper and is presented in details in a separate paper; see \cite{GazorMokhtari} and also \cite{CushmanSanders,CushSandCont}. We also would like to acknowledge Professor James Murdock's generous help, discussions, and remarks through the first author in his 2011 summer visit and numerous email communications. The following lemma represents the second level normal form of volume-preserving vector fields of Hopf-zero singularity, where a \(\Sl\)-style normal form is applied.

\begin{lem} \label{2ndLevel} There exists an invertible transformation transforming \(v^{(1)}\) into
its second level normal form
\be
v^{(2)}= \Theta^0_0+\alpha_0 F^{-1}_0+ \sum^\infty_{k=1} \alpha_k^{(2)} F^k_k+\sum^\infty_{k=1} \beta_k^{(2)}\Theta^k_k,
\ee where \(\alpha_0= a^{-1}_0.\)
\end{lem}
\bpr Define a grading function by \(\delta(F^l_k):=\delta(\Theta^l_k):= k.\) Then, the result is deduced from $[ F^{-1}_0, F^l_k]
= 2(l+1)F^{l-1}_{k}$ and \([F^{-1}_0, \Theta^l_k]= 2l\Theta^{l-1}_k.\)
\epr
Note that by a linear change of state variable, we can remove \(\Theta^0_0\) from the system; see \cite[Theorem 4.1]{GazorMokhtariEul}
and \cite[Lemma 5.3.6]{MurdBook}.
Assume that there exist \({\alpha_k^{(2)}}\neq0\) (for \(k\geq 1)\) and denote
\be\label{PQDef} p:=\min\left\{k\,|\,
\alpha^{(2)}_k\neq0, k\geq 1\right\}.
\ee
Let
\bes \alpha_0:=\frac{1}{2}, \alpha_p:= \alpha_p^{(2)} \hbox{ and } {\X_p}:=\frac{1}{2}F^{-1}_{0}+\alpha_{p}F^{p}_{p}.\ees
Define a grading structure by
\bes
\delta\left(F^l_k\right):= p(k-l)+k\;\;
\hbox{ and }\;\;
\delta\left(\Theta^l_k\right):= p(k-l)+k+p+1.
\ees

\begin{rem}\label{arbitA}
For any \(a, b\in \f, a, b>0, \) through linear changes of variables
\bes t:=\frac{1}{b\alpha_0}\left(\frac{ab\alpha_0\, {\rm sign}(\alpha_0\alpha_p)}{\alpha_p}\right)^\frac{1}{p+1}\tau, \
x:= \left(\frac{ab\alpha_0\, {\rm sign}(\alpha_0\alpha_p)}{\alpha_p}\right)^\frac{1}{p+1}X, \ y:= Y, \ \hbox{ and } \ z:= Z,\ees  we can transform \(\X_p\) into
\(\X_p:= bF^{-1}_0+ a\,{\rm sign}(\alpha_p\alpha_0)F^p_p.\) Consequently, the coefficients of \(F^{-1}_0\) and \(F^p_p\) can be arbitrarily
chosen in the normal form
computation. Thus without the loss of generality, we can choose \(\alpha_0:=\frac{1}{2}.\) Hence, we can change the coefficient
\(\alpha_p\) into a non-algebraic number. We shall use this in Theorem \ref{FinalThm} in order to simplify the system as desired, and then again we change \(\alpha_p\) into \(\pm 1.\) As a theoretical result, this may not violate the principles
of normal form theory, but this needs more attention when it is implemented in a computer program.
In fact irrational numbers are treated like rational numbers in computers because of computers' round off errors.
This, however, does not hamper our results. Indeed for implementation of the results on any
computer, one needs to truncate the system up to a certain degree and thus, one only needs to choose \(\alpha_p\)
to be distanced from the roots of a finite number of polynomials. For one of the three cases (considered by Baider and Sanders \cite{BaidSand91}) of Bogdanov--Takens singularity when certain ratio of coefficients is non-algebraic, the simplest normal form is known; see \cite{KokubuWang}. This case is still an open problem when the ratio is algebraic. Hence, this technique is very useful wherever it is applicable. Indeed, we believe that it can be applied to other problems; \eg see \cite{GazorMokhtari}.
\end{rem}

The following lemma introduces the transformation needed for elimination of the term \(\Theta^m_n\) when \(m<n.\)

\begin{lem}\label{convert}
For natural numbers \(m\) and \(n,\) there is a \(\delta\)-homogenous polynomial state solution \(\Y^m_n\) such that
\be\nonumber
\Theta^m_n+ \left[\Y^m_n,\X_p\right] = \frac{{\alpha_p}^{n-m}\big(p(m-n+1)-n\big)\big(m+1+(m-n+1)(p+1)\big)^{n-m-1}_{2p+2}}{(-2)^{n-m}(m+1)^{n-m}_{p+1}}
\,\Theta^{pn-pm+n}_{pn-pm+n}.
\ee
\end{lem}
\bpr The proof follows by choosing
\bes
\Y^m_n:=\sum _{i=0}^{n-m-1} \frac{(-1)^{i+1}{\alpha_{{p}}}^{i}\big(m+1+(m-n+1)(p+1)\big)^i_{2p+2}}{{2}^{i+1}(m+1)(m+p+2)^i_{p+1}}
\Theta^{m+ip+i+1}_{n+ip}.
\ees
\epr

\begin{lem}\label{p+1th}
The \((p+1)\)-th level normal form of \(v\) associated with Equation \eqref{Eq1} is
\be v^{(p+1)}=\Theta^0_0+\frac{1}{2} F^{-1}_0+\alpha_{p}F^{p}_{p}+ \sum\alpha_k^{(p+1)} F^k_k+\sum\beta_k^{(p+1)}\Theta^k_k,
\ee where the first summation is over \(k\neq p-1 \mod{2(p+1)}\) and the second summation is over \(k\neq -1 \mod{2(p+1)}.\)

\end{lem}
\bpr Since \(F^k_k, \Theta^k_k\in \ker \ad_{F^1_0},\) we follow Baider and Sanders \cite{BaidSand91} and define
\(\mathscr{G}:=\ad(F^1_0)\circ\ad({\X_p}).\) Then,
\bas
\mathscr{G} (F^{l}_{k})&=&4(l+1)(l-k-2)F^{l}_{k}+2\alpha_{p}(k-l+1)
\big(p-l+(k-l)(p+1)\big) F^{p+l+1}_{p+k},
\\
\mathscr{G}(\Theta^{m}_{n})&=& 4m(m-n-1)\Theta^{m}_{n}+2\alpha_{p}(m-n)\big(m(p+2)-n(p+1)\big)
\Theta^{p+m+1}_{p+n}.
\eas
Let
\bes \left\{F^l_k, \Theta^m_n \,|\,-1\leq l\leq k, 0\leq m\leq n\right\}\ees
be an ordered basis for \(\LST,\) where its ordering is partially defined by \(F^l_k\prec F^m_n\) and
\(\Theta^l_k\prec\Theta^m_n\)
if \(k<n.\)
Now the matrix representation of \(\mathscr{G}\) is lower
triangular. Thus, for any natural number \(k\) there exist \(\delta\)-homogenous polynomial vector fields
\be
\left\{
  \begin{array}{ll}
\mathcal{F}^{-1}_k&:=\sum_{m=0}^{2k}\frac{{\alpha_p}^{m}(2k+1)^{m-1}_{-2}}{{2}^{m}(m)!} F^{mp+m-1}_{2k+mp-1}, \qquad
\delta\left(\mathcal{F}^{-1}_k\right) \equiv -1 \mod{2(p+1)},
\\
\XS_p^k&:= \sum^k_{m=0} \binom{k}{m} \frac{{\alpha_p}^{k-m}}{2^m} F^{m(p+1)-1}_{2k+mp+2}, \quad \qquad
\delta\left(\XS^k_p\right)\equiv p \mod{2(p+1)},
\\
\label{KerTheta}\mathcal{T}^0_{k, p}&:=\sum_{m=0}^{k}{\frac{{\alpha_{{p}}}^{m}(k)_2^m}{ {2}^{m}m! }}
\Theta^{m(p+1)}_{k+mp}, \qquad \qquad\;\;\;
\delta\left(\mathcal{T}^0_{k, p}\right)\equiv 0 \mod{ p+1},
\end{array}
\right.
\ee
so that
\be
\ker(\mathscr{G})= \Span \left\{\mathcal{F}^{-1}_k, \XS_p^{k}, \mathcal{T}^0_k\,|\, k\in \N\right\}.
\ee
On the other hand,
\bas
\left[\X_p, \XS^k_p\right]&=&0,
\\
\left[{\X_p},\mathcal{T}^0_{k,p}\right]&=&\frac{k(p+1){\alpha_{p}}^{k}(k)^k_{-2}}{ {2}^{k}k! }\Theta^{k(p+1)+p}_{k(p+1)+p},
\\
\left[\X_p, \mathcal{F}^{-1}_k\right]&=&\frac{(2k-1)(p+1){\alpha_{{p}}}^{2k+1}(2k-1)^{2k-1}_{-2}}{{2}^{2k}(2k)!} F^{2k(p+1)+p-1}_{2k(p+1)+p-1}.
\eas

\noindent Therefore, \(F^m_m, \Theta^n_n\in \IM d^{m, p+1}\) for any \(m\equiv p-1 \mod{2(p+1)}\) and \(n\equiv -1 \mod{2(p+1)},\) where \(m=n+p+1.\)
This completes the proof.
\epr

The following theorem presents the simplest normal form for the system \eqref{Eq1}.

\begin{thm}\label{FinalThm}
There exist invertible transformations (including linear time rescaling)
sending \(v\) given by Equation \eqref{Eq1} into the \((p+q+2)\)-th level
normal form system
\ba\label{Inft}
\left\{
  \begin{array}{ll}
\dot{x}&=(y^2+z^2)\pm  x^{p+1}+x^{p+1}\sum^\infty_{k=1} \alpha_{k+p} x^{k},\\
\dot{y}&=zg(x)\mp\frac{(p+1)}{2} x^{p}y-x^py\sum^\infty_{k=1}\frac{(k+p+1)\alpha_{k+p}}{2}x^{k}, \\
\dot{z}&=-yg(x)\mp\frac{(p+1)}{2}x^{p} z-x^pz\sum^\infty_{k=1}\frac{(k+p+1)\alpha_{k+p}}{2}x^{k},
  \end{array}
\right.
\ea for \(p, q\geq 1,\) and
\bes
g(x):=1+ x^q\sum^\infty_{k=0} \beta_{q+k}x^k,
\ees where \(\alpha_{k+p}=0\) for \(k\equiv -1 \mod{2(p+1)}.\) Furthermore, \(\beta_{q+k}=0\) for \(k\equiv p \mod{2(p+1)}\) and
for \(k\equiv -(q+1) \mod{2(p+1)}.\) In addition, the \((p+q+2)\)-th level normal form system \eqref{Inft} is the infinite level normal form. Let
\bas
f(x, y, z)&:=&\left(y^2+z^2\right)\left( \frac{1}{2} \left(y^2+z^2\right)\pm x^{p+1}+x^{p+1}\sum^\infty_{k=1}\alpha_{k+p} x^k\right).
\eas and \(\F\) denote for the algebra of first integrals of \(v^{(\infty)}.\) Then, \(\F= \langle f\rangle\)
and the symmetry group of \(v^{(\infty)}\) is generated by \(\F\Theta^0_0.\)
\end{thm}
\bpr Assume that \(\beta_i^{(p+1)}\neq 0\) for some \(i\geq 1.\) Then, define
\bes
q:=\min\left \{i\,|\, \beta^{(p+1)}_i\neq0, i\geq 1\right\}.
\ees For any \(k> 0\) we have
\bes \XS^{k+1}_p=\sum_{m=0}^{k+1}{\alpha_p}^m {{k+1}\choose{m}} F^{m(p+1)-1}_{2k+mp},\ees
and by
Lemma \ref{convert} there exists a state solution \(Y\) such that
\ba
&&\left[ \XS^{k+1}_p,\Theta^q_q\right]+\left[Y, \X_p\right]
\\\nonumber&=&
\sum_{m=0}^{k+1}
\frac{q{{k+1}\choose{m}}(2k+2-m){\alpha_p}^{2k-m+1}\big(q+2(m-k)(p+1)\big)^{2k-m+1}_{p+1}}{(-1)^{m+1}{2}^{2k-m+1}\big(q+j(p+1)\big)^{2k-m+1}_{p+1}}\,
\Theta^{2k(p+1)+p+q}_{2k(p+1)+p+q}.
\ea
By Remark \ref{arbitA}, without any loss of generality, we may assume that \(\alpha_p\) is not an algebraic number and thus,
\bes \Theta^m_m\in \IM d^{m+p+1, p+q+2}, \hbox{ where } m\equiv p+q\mod{2(p+1)}.\ees For \(k=2l,\) we have
\bas
\label{KerTheta}\mathcal{T}^0_{2l,p}
&=&\left(y^2+z^2\right)^l\left(y^2+z^2+\alpha_px^{p+1}\right)^{l}\left(z\ddy-y\ddz\right)\in \ker {\rm ad}_{\X_p}.
\eas Then, \(\mathcal{T}^0_{2l,p}\) is extended to a symmetry for \(\X_p+\sum^n_{k=p+1}\alpha_{k} F^{k}_{k},\) \(n\in \N,\) \ie
\bas
\label{KerTheta}
\mathcal{T}^0_{2l,n}&:=&
\left(y^2+z^2\right)^l\left(y^2+z^2+\sum^n_{k=p}\alpha_{k}x^{k+1}\right)^l\left(z\ddy-y\ddz\right)\in \ker {\rm ad}_{\X_p+\sum^n_{k=p}\alpha_{k} F^{k}_{k}}.
\eas Since \(\mathcal{T}^0_{2l,n}\) is convergent to \(f^l\Theta^0_0\) with respect to the filtration topology, the proof is complete. \epr

\section{Parametric normal forms }\label{PNF}

In this section, we deal with the parametric normal form of a multiple-parametric perturbation of the system \eqref{Eq1}.
Roughly speaking, this section provides an infinite level parametric normal form for the miniversal unfolding for the system \eqref{Eq1}; also see \cite{GaoPNF,GazorYuSpec,GazorYuFormal,Murd04,Murd08} and Remark \ref{shift}. Since nonlinear time rescaling destroys the symmetry of the system
(\ie volume-preserving), we do not use nonlinear time rescaling.
However, parametric time rescaling is permitted when it does not depend on the state variables. Consider a parametric vector field
\ba\label{PEq}
w(x, y, z, \mu):= \sum a^l_{k,\m}F^l_k\mu^\m+ \sum b^i_{j,\m} \Theta^i_j\mu^\m,
\ea where
\ba\nonumber
&k\geq -1, k\geq l\geq -1, j\geq i\geq 0,&\\\nonumber
&\m=\left(n_1, n_2, \ldots, n_m\right), n_j\in \N\cup\{0\} \hbox{ for } j= 1, 2, \ldots, n_m, &\\
&\mu:=\left(\mu_1, \mu_2, \ldots, \mu_m\right), \mu^\m:= \mu_1^{n_1}\mu_2^{n_2}\ldots \mu_m^{m_m}, \sum^m_{j=1} n_j\geq 0,&\\\nonumber
&b^0_{0,\0}\neq 0, a^{-1}_{0,\0}\neq 0, a^{-1}_{-1,\0}=0,\hbox{ and } a^0_{0,\0}=0.&
\ea
We call any parametric vector field \(w\) given in Equation \eqref{PEq}, a multi-parametric
deformation for \(v\) when \(v=w(x, y, z, \0).\) By a similar argument to the one used in the proof of Lemma \ref{2ndLevel}
and a parametric time rescaling, any multi-parametric deformation \(w\) for \(v\) associated with Equation \eqref{Eq1}
can be transformed to the second level extended partial parametric normal form
\bas
w^{(2)}=\sum \beta_{0,\m}^{(2)}\Theta^0_0\mu^\m+\sum \alpha^{(2)}_{0,\m}F^{-1}_0\mu^\m+ \sum \alpha_{k,\m}^{(2)} F^k_k\mu^\m
+\sum \beta_{j,\m}^{(2)}\Theta^j_j\mu^\m,
\eas where \bes \beta_{0,\0}^{(2)}=1, \alpha^{(2)}_{0,\m}=\frac{1}{2}, k\geq -1, \hbox{ and } j\geq 1; \ees see \cite{GazorYuSpec} for more details. We consider a parametric change of state variable \bes [x,y,z]:=exp\left(h(\mu)\left(tY\ddz-tZ\ddy\right)\right)[X,Y,Z];\ees this is similar to the argument used in \cite[Theorem 4.1]{GazorMokhtariEul} and \cite[Lemma 5.3.6]{MurdBook}.
Here, \bes h(\mu):= 1+ \sum \beta^{(2)}_{0, \m}\mu^\m, \hbox{ and } [x, y, z]:= x\ddx+y\ddy+z\ddz\ees denotes the new variables and \([X,Y,Z]\) stands for the old variables. Once all linear terms are omitted, using a parametric time rescaling, we can transform \(w^{(2)}\) into
\be\label{PNF2NL}
\tilde{w}^{(2)}= \frac{1}{2}F^{-1}_0+ \sum \tilde{\alpha}_{k,\m}^{(2)} F^k_k\mu^\m+\sum \tilde{\beta}_{j,\m}^{(2)}\Theta^j_j\mu^\m,
\ee where \(k\geq -1\) and \(j\geq 1.\) Let
\bes \tilde{\alpha}_{k,\0}^{(2)}\neq0 \hbox{ for some } k\geq 1, \hbox{ and } \tilde{\beta}^{(2)}_{j,\0}\neq 0 \hbox{ for some } j\geq 1.\ees
Denote
\bes p:=\min\left\{k\,|\, \tilde{\alpha}^{(2)}_{k,\0}\neq0, k\geq 1\right\} \hbox{ and }\; q:=\min\left\{j\,|\,
\tilde{\beta}^{(2)}_{j,\0}\neq0, j\geq 1\right\}.
\ees We define the grading function \(\delta\) by
\bas
\delta\left(F^l_k\mu^\m\right)&:=&p(k-l)+k+(p+q+3)|\m|,
\\
\delta\left(\Theta^l_k\mu^\m\right)&:=& p(k-l)+k+p+1+ (p+q+3)|\m|.
\eas

\begin{lem}\label{2ndLevelP}  There exist invertible changes of variables that they transform \(\tilde{w}^{(2)}\) given by
Equation \eqref{PNF2NL} into
the \((p+q+2)\)-th level extended partial parametric normal form
\be\label{pq2LvlP}
w^{(p+q+2)}:= \frac{1}{2} F^{-1}_0+ \sum \alpha_{k,\m} F^k_k\mu^\m+\sum \beta_{j,\m}\Theta^j_j\mu^\m,
\ee where \(k\geq -1, j\geq 1, \alpha_{k,\0}=0\) for all \(k<p,\) \(\alpha_{p,\0}=\pm1,\) \bes \alpha_{k,\m}=0\ \hbox{ for } \ k\equiv p-1 \mod{2(p+1)},\ees
and \(\beta_{k,\m}=0\) for both \bes k\equiv -1 \mod{2(p+1)}\ \hbox{ and } \ k\equiv -(q+1) \mod{2(p+1)}.\ees
\end{lem}

\bpr Note that the number \(q\) must be updated in the (\(p+1\))-th level normal form by
\bes q:=\min\left\{j\,|\, \beta^{(p+1)}_{j,\0}\neq0, j\geq 1\right\}.\ees
The proof is straightforward by similar arguments in the proofs of Lemma \ref{p+1th} and Theorem \ref{FinalThm}.
\epr

By omitting terms of degree (standard degree of polynomials) higher than or equal to \(n+1\) of a vector field \(v,\) we obtain its \(n\)-degree
truncated (\(n\)-jet) vector field and denote it by \(J^n(v).\)
For any natural number \(N,\) let \bes k:= \left\lfloor \frac{N}{2(p+1)}\right\rfloor, \ \ l:= N-2k(p+1),\ees
and \bes r:= k(2p+1)+l+2,\ees where \(\lfloor a\rfloor\) denotes the integer part of the real number \(a\).
Denote \(n_j\) (\(j=2, \ldots, r\)) for all natural numbers in which
\be\label{nj}
-1\leq n_j\leq N-1 \hbox{ and, for all } n_j>p  \hbox{ we have } n_j\neq p-1 \mod{2(p+1)}.
\ee
Denote \(m_{k}\) (\(k=1,\ldots s\)) for all natural numbers
\be\label{mj} 1\leq m_k\leq N, \ \hbox{ where } m_k\neq -1 \mod{2(p+1)}\hbox{ and } m_k\neq p+q \mod{2(p+1)}.\ee
Obviously,
\bes n_1:= -1, n_2:=0, n_3:=1 \hbox{ and } m_1:=1.\ees
Now for the parametric normal form \(w^{(p+q+2)}\) given in Lemma \ref{2ndLevelP}, we denote
\bes
A_N:= (a_{ij}), \hbox{ where } a_{ij}= \alpha^{(p+1)}_{n_j,\mu_i} \hbox{ for } j=1, \ldots, r, \hbox{ and } a_{ij}= \beta_{m_{j-r},\mu_i}^{(p+1)}
\hbox{ for } j=r+1, \ldots, r+s.
\ees
The truncated normal form \(J^{N+1}(w^{(p+q+2)})\) is called a {\em non-degenerate perturbation} when
\be\label{Rank}\rank(A_N) = r+s.\ee

\begin{rem}\label{shift} When \(N=p,\) a transformation of the form \(x=X+\epsilon\) (for an appropriate formal function \(\epsilon=\epsilon(\mu)\)) can be applied to Equation \eqref{pq2LvlP} in order to simplify \(\sum \alpha_{p-1, \m}F^{p-1}_{p-1}\mu^\m.\)  Then, this reduces the number of parameters left at the truncated parametric normal form. Thereby, one needs a weaker rank condition than the rank condition \eqref{Rank} in order that a \(p\)-degree truncated parametric normal form system would be parametric generic. Indeed, the sequence \(n_j\) skips \(p-1\) and we have \(r:= k(2p+1)+l+1.\) Then, \(\rank(A_p):= r+s.\) These types of transformations was extensively discussed by Murdock and Malonza \cite{Murd08}. However, we shall not use them in this paper.
\end{rem}

The condition \eqref{Rank} guarantees that a reparametrization sends the parametric vector field \(J^{N+1}(w^{(p+q+2)})\) into the \((N+1)\)-degree truncated
infinite level parametric normal form
\ba
J^{N+1}(w^{(\infty)}):= \frac{1}{2} F^{-1}_0\pm F^p_p+ \sum^r_{i=1} F^{n_i}_{n_i}\mu_i+ \sum^N_{n=p+1} \alpha_{n} F^{n}_{n}
+ \sum^N_{k=q} \beta_{k} \Theta^k_k+ \sum^{r+s}_{i=r+1} \Theta^{m_{i-r}}_{m_{i-r}}\mu_i.
\ea
This motivates the following theorem.
\begin{thm}\label{InftPNF}
There exist formal invertible changes of state variables and parametric
time rescalings that send any non-degenerate deformation of the system \eqref{Eq1} into the \((p+q+2)\)-th level parametric normal form \(w^{(p+q+2)}.\) Furthermore, a reparametrization sends \(J^{N+1}(w^{(p+q+2)})\) into the \((N+1)\)-degree truncated infinite level parametric normal form
\ba\label{PInft}
\left\{
  \begin{array}{ll}
\dot{x}&= (y^2+z^2)\pm  x^{p+1}+\sum^r_{i=1} x^{n_i+1}\mu_i+x^{p+1}\sum^{N-p}_{k=1} \alpha_{k+p} x^{k},\\
\dot{y}&=zg(x)\mp\frac{(p+1)}{2} x^{p}y-\sum^r_{i=1} \frac{n_i+1}{2}x^{n_i}y\mu_i-x^py\sum^{N-p}_{k=1}\frac{(k+p+1)\alpha_{k+p}}{2}x^{k}, \\
\dot{z}&=-y g(x)\mp\frac{(p+1)}{2}x^{p} z-\sum^r_{i=1} \frac{n_i+1}{2}x^{n_i}z\mu_i-x^pz\sum^{N-p}_{k=1}\frac{(k+p+1)\alpha_{k+p}}{2}x^{k},
  \end{array}
\right. \ea where
\ba\label{gPar}
g(x)&:=&1+ \sum^{r+s}_{i=r+1} x^{m_{i-r}}\mu_i+ x^q\sum^{N+1-q}_{k=0} \beta_{q+k}x^k, \quad p, q\geq 1,
\ea the coefficients \(\alpha_{k}=0\) for \( k\equiv p-1 \mod{2(p+1)},\) and
\bas
\beta_{l}=0& \hbox{ for }& l\equiv -1 \mod{2(p+1)} \ \hbox{ and } \ l\equiv p+q \mod{2(p+1)}.\eas
Let \(\F\) be the algebra of first integral for \(w^{(\infty)}.\) Then,
\ba\label{P1integ}
\F&=&\left\langle \left(y^2+z^2\right)\left( \frac{1}{2} \left(y^2+z^2\right)\pm x^{p+1}+\sum^r_{i=1} x^{n_i+1}\mu_i+x^{p+1}\sum^{N-p}_{k=1}\alpha_{k+p} x^k\right)\right\rangle.
\ea The parametric symmetry group of \(w^{(\infty)}\) is generated by \(\F\Theta^0_0.\)
\end{thm}

\bpr
The proof is complete by applying a linear change of state variables
\bes [x,y,z]:=\exp\left(\left(Y\ddz-Z\ddy\right)t\right)[X,Y,Z]\ees
to \(w^{(p+q+2)}\),
where \(X,Y,Z\) are the new variables and \(x,y,z\) show the old variables. This returns the omitted linear part \(\Theta^0_0\) back into the system.
\epr

\section{Examples }\label{Examples}

In this section, we first derive several necessary relations between coefficients. These are enough so that the cubic-degree truncated classical normal form system would belong to the space \(\LST.\) We also obtain some useful formulas for the infinite level normal forms of a generic cubic classical normal form system \(v\in \LST.\)
Then, we apply these relations to obtain a one-parameter family of modified R\"{o}ssler and a one-parameter class of
generalized Kuramoto--Sivashinsky equations.
Then, we apply our results to these systems. The system of Hopf-zero singularity can have very rich dynamics.
Some major contributions have been made in the literature; see \eg
\cite{AlgabaHopfElect,LlibreHopfZeroBif,DumortierHopfZ,Buzzi,HarlimLangf,LangfTori,LangfHopfSteady,LangfHopfHyst}. However, there are many degenerate or symmetric cases that their dynamics have not yet been investigated.

Consider a differential system governed by
\be\label{PreNF}
\left(\begin{array}{ccc}\dot{x} \\\dot{y} \\\dot{z} \\\end{array}\right)=\left(\begin{array}{ccc}0& 0& 0 \\0&0 &1 \\0&-1 & 0
\\\end{array}\right)\left(\begin{array}{c}
x\\y\\z\end{array}
\right)+\sum_{2\leq i+j+k}\left(\begin{array}{ccc}a_{ijk}\\b_{ijk}
\\c_{ijk}\end{array}\right)x^iy^jz^k.
\ee
\noindent Let \(v^{(1)}\) denote the cubic-truncated classical normal form of the system \eqref{PreNF}. Now the following three relations
are the necessary conditions for \(v^{(1)}\) to be in \(\LST\):
\ba
0&=&a_{002}+b_{{011}}+c_{{101}},\\\nonumber
0&=&2a_{003}+3b_{{012}}+3 c_{{102}}-b_{{002}}(3b_{{110}}-4a_{{101}}+6c_{{200}})-c_{{002}}(4a_{{011}}-6b_{{020}}
-3c_{{110}})
,\\\nonumber
0&=& b_{{020}}+16a_{021}+2b_{{210}}+2c_{{200}}+2c_{{020}}+8 a_{{102}}-16b_{{020}} a_{{101}}+16c_{{020}} a_{{011}}
\\&&\nonumber
-7c_{{011}}(a_{{200}}-a_{{020}})- 16a_{{110}}(c_{{101}}+7b_{{011}})-b_{{101}}(7a_{{200}}-15a_{{020}}) \\\nonumber
&& - 2b_{{200}}(b_{{110}}+2c_{{200}}+2b_{{020}}+8a_{{101}})+ 2c_{{200}}(c_{{110}}+8a_{{110}}+8a_{{011}}).
\ea Then,
\be\label{CubicNF}
v^{(1)}=\Theta^0_0+a^{-1}_0 F^{-1}_0+a^1_1F^1_1+a^0_1F^0_1+a^2_2F^2_2+
b^1_1\Theta^1_1+b^0_2\Theta^0_1+b^2_2\Theta^2_2,
\ee where
\ba\nonumber
a^{-1}_0 &=&\frac{1}{2}(a_{0 2 0}+a_{2 0 0}),
\\\nonumber
a^1_1&=&\frac{1}{2}(b_{{011}}+c_{{101}}),
\\\nonumber
a^0_1&=&\frac{1}{32}\big(b_{{020}}+2b_{{210}}+2c_{{200}}+2c_{{020}}
+c_{{011}}(a_{{200}}-a_{{020}})+a_{{110}}(b_{{011}}
-c_{{101}})\big)
\\\nonumber
&&+\frac{1}{16}c_{{200}}(c_{110}-2b_{{200}})
+\frac{1}{32}b_{{101}}(a_{{200}}-a_{{020}})
-\frac{1}{16}b_{{200}}(b_{{110}}+2b_{{020}}),
\\\nonumber
a^2_2&=&\frac{1}{2}\big(c_{{102}}+b_{{012}}+c_{{002}}(c_{{110}}-2a_{{011}}+2b_{{020}})+b_{{002}}(2a_{{101}}-2c_{{200}}
-b_{{110}})\big).
\ea These equations are derived from the fundamentally useful formulas given in \cite{AlgabaHopfZ}, where explicit formulas for \(b^1_1,b^0_2,\)
and \(b^2_2\) are also given.

\begin{prop}\label{Propcoef}
Consider the vector field
\bas
v:= \Theta^0_0+a^{-1}_0{F^{-1}_0}+a^0_1F^0_1+a^1_1F^1_1+a^2_2F^{2}_2+b^1_1\Theta^1_1+b^0_1\Theta^0_1+
b^2_2\Theta^2_2,
\eas where \(a^{-1}_0, a^1_1\neq 0.\) Then, the quartic truncated infinite level normal form is governed by
 \bas
{v}^{\infty}:=\Theta^0_0+ \frac{1}{2} F^{-1}_0\pm F^{1}_1
\pm \frac{a^{-1}_0 a^2_2}{\sqrt{2}a^1_1|a^{-1}_0a^1_1|^{\frac{1}{2}}}
F^{2}_2+ \frac{a^0_1a^2_2}{8{a^1_1}^2}F^3_3\pm\frac{b^1_1}{a^1_1}\Theta^1_1
\pm \frac{4b^2_2 a^{-1}_0+a^0_1b^1_1}{4\sqrt{2}a^1_1|a^{-1}_0a^1_1|^{\frac{1}{2}}}\Theta^2_2,
\eas where \(\pm 1\) represents \({\rm sign}(a^{-1}_0a^1_1).\)
\end{prop}
\bpr
By state change of variables, we can transform the system into
 \bas
\tilde{v}:=a^{-1}_0{F^{-1}_0}+a^1_1{F^{1}_1}+a^2_2{F^{2}_2}+\frac{a^0_1a^2_2}{4a^{-1}_0}F^3_3+b^1_1\Theta^1_1+\frac{4b^2_2
a^{-1}_0+a^0_1b^1_1}{4a^{-1}_0}\Theta^2_2.
 \eas The linear changes of variables \bes t:= \left(\frac{2 \,{\rm sign}(a^{-1}_0a^1_1)}{a^{-1}_0a^1_1}\right)^\frac{1}{2}\tau, \
x:= \left(\frac{a^{-1}_0\, {\rm sign}(a^{-1}_0a^1_1)}{2a^1_1}\right)^\frac{1}{2}X, \ y:= Y, \ \hbox{ and } \ z:= Z\ees transform this system into the desired form. Next, a second linear changes of variables is needed to add \(\Theta^0_0\) back into the system.
\epr

Now we apply the above formulas to two examples.

\begin{exm}\label{RosselerExample}
Consider a modified R\"{o}ssler equation governed by
\bas
\dot{x}&=&-y-z+dy^2,\\
\dot{y}&=&x+ay+ez^3,\\
\dot{z}&=&bx-cz+xz+rz^3.
\eas
\noindent The parameter values \(d=e=r=0\) gives rise to the well-known R\"{o}ssler equation and for two sets of parameter values, it
has Hopf-zero singularity at origin; one is for parameter values \(a=c, b=1, a^2<2,\) while the other has simple dynamics; see \cite{AlgabaHopfZ}.
Hence, we choose
\bes
a:=c, b:=1, d:=a^2-1, \hbox{ and } e:=-2{a}^{3}+{\frac {15}{17}}r{a}^{2}+{\frac {5}{17}}a-{\frac {13}{17}}r,
\ees where
\bes r:=\frac {510{a}^{7}-891{a}^{5}-170{a}^{4}-1316{a}^{3}+510{a}^{2}+1058a-340}{15({a}^{2}-2)(15{a}^{4}-2{a}^{2}-96)} \hbox{ and } -\sqrt{2}< a<\sqrt{2}.\ees

By a linear change of coordinates we may transform the linear part into the Jordan canonical form and then the system is given by
\ba\nonumber
\dot{x}&:=&\frac{-1}{2(2-{a}^{2})^\frac{3}{2}}\bigg(2(2-{a}^{2})^\frac{3}{2}y+({a}^{2}-2d){x}^{2}
-\sqrt {2-{a}^{2}}({a}^{2}+4d)xy-a({a}^{2}-4d+2)xz
\\\nonumber&&+2d({a}^{2}-2){y}^{2}+2a\sqrt {2-{a}^{2}}(2d+1)yz-2{a}^{2}(d-1){z}^{2}
-3a(r-e)\sqrt{2-{a}^{2}}y({x}+{a}{z})^2
\\\nonumber&&+a({a}^{2}-2)(r-e)y^2(
\sqrt {2-{a}^{2}}{y}-3x+3az)
+a(r-e)({x}-{a}z)^3
\bigg),
\\\label{RosselerEqn}
\dot{y}&:=&\frac{1}{2({a}^{2}-2)}\Big(2({a}^{2}-2)z-a{x}^{2}+(2+{a}^{2})xz
-2a{z}^{2}-3\sqrt {2-{a}^{2}}(r+e)y({x}-z)^2
\\\nonumber&&
+({a}^{2}-2)(r+e)y^2(\sqrt {2-{a}^{2}}{y}+3az-3x)
+(r+e)(x-a{z})^{3}
+\sqrt {2-{a}^{2}}y(ax-2z)
\Big),
\\\nonumber
\dot{z}&:=&\frac{-1}{(2-{a}^{2})^\frac{3}{2}}\bigg(-a(d-1){x}^{2}+\sqrt {2-{a}^{2}}y\big(2(d{a}^{2}+1)z-a(2d+1)x\big)
+(2d{a}^{2}-{a}^{2}-2)xz
\\\nonumber&&+ad({a}^{2}-2){y}^{2}
-a(d{a}^{2}-2){z}^{2}+({a}^{2}-2)(r-e)y^2\big(\sqrt {2-{a}^{2}}y+3az-3x\big)
\\\nonumber&&+(r-e)({x}-az)^{3}-3(r-e)\sqrt{2-{a}^{2}}y(x-az)^2\bigg).
\ea
Then, the vector field
\be\label{CubicNFRos}
v^{(1)}=\Theta^0_0+a^{-1}_0F^{-1}_0+a^1_1F^1_1+a^0_1F^0_1+a^2_2F^2_2+b^0_1\Theta^0_1 +
b^1_1\Theta^1_1+b^2_2\Theta^2_2
\ee is the cubic-truncated classical normal form, where
\bes a^{-1}_0:=-\frac {a}{\sqrt{2-a^2}}, a^1_1:= \frac {a}{2\sqrt {2-a^2}},
a^2_2:=\frac{15{a}^{7}-24{a}^{5}-5{a}^{4}-39{a}^{3}+15{a}^{2}+17a-10}{2.5 (2-{a}^{2})^{3/2}(15{a}^{4}-2{a}^{2}-96)},
\ees
\ba\label{F01}
a^0_1&:=& \frac{15{a}^{9}+577{a}^{7}+60{a}^{6}-1788{a}^{5}-420{a}^{4}
+504{a}^{3}+840{a}^{2}+240a-480}{40(15{a}^{4}-2{a}^{2}-96)(2-{a}^{2})^{3/2}},
\\\nonumber
b^0_1&:=& \frac{4185{a}^{8}-7671{a}^{6}-1020{a}^{5}-12331{a}^{4
}+3060{a}^{3}+19938{a}^{2}-2040a-15840}{-240(a^2-2)^{2}(15{a}^{4}-2{a}^{2}-96)},
\\\nonumber
b^2_2&:=& \frac {17235{a}^{8}-28066{
a}^{6}-2720{a}^{5}-72356{a}^{4}+8160{a}^{3}+96888{a}^{2}-5440a-32640}{-320(a^{2}-2)^{2}(15{a}^{4}-2{a}^{2}-96)},
\\\nonumber
b^1_1&:=&\frac{3{a}^{2}-1}{4{a}^{2}-2}.
\ea Thereby, the infinite level normal form for any arbitrary \(a\neq 0\) (\(a^2<2\))
is obtained through Proposition \ref{Propcoef}.
\end{exm}

\begin{exm}\label{KSExample}
 Consider a generalized Kuramoto--Sivashinsky equation given by
\be\label{SK}
h_t+\ddx G(h, h_x, h_{xx})=0,
\ee where \bes G\left(h, h_x, h_{xx}\right):=2h^2+ h_x+a{h_x}^2+b{h_{xx}}^2+c{h_x}^2h_{xx}+d h_{xx}h^2+h_{xxx}. \ees
A common approach is to consider
traveling wave solutions \(h(x, t)=u(x+\alpha t).\) Then, this gives rise to an ordinary differential equation
\bes
\alpha u'+4uu'+u''+ 2a u'u''+2d uu'u''+2c u'{u''}^{2} + 2b u''u^{(3)}+ d u^2u^{(3)}+c\left(u'\right)^2u^{(3)}+u^{(4)}=0.
\ees Integrating this equation (assuming that \(u=0\) is a solution for all \(\alpha\)), we have
\be\label{SKODE}
\alpha u+2u^2+u'+a\left(u'\right)^2+b\left(u''\right)^2+c\left(u'\right)^2u''+d u''u^2+u^{(3)} =0.
\ee

\noindent Let \(\alpha :=0\) (this is to get a Hopf-zero singularity) and \(v:=(v_1, v_2, v_3)= ( u,u',u'').\)
Then, using a linear change of complex variables, we may send the system to
\be\label{KS}
\left\{
  \begin{array}{ll}
\dot{x}&=-(x+y+z)^{2}\big(2+d(y+z)\big)-(y-z)^{2}\big(a-c(y+z)\big)+b(y+z)^{2},\\
\dot{y}&=iy+(x+y+z)^{2}\big(1+\frac{d}{2}(y+z)\big)+\frac{1}{2}(y-z)^{2}\big(a-c(y+z)\big)-\frac{b}{2}(y+z)^{2},\\
\dot{z}&=-iz+(x+y+z)^{2}\big(1+\frac{d}{2}(y+z)\big)+ \frac{1}{2}(y-z)^{2}\big(a-c(y+z)\big)-\frac{b}{2}(y+z)^{2}.
  \end{array}
\right.
\ee Here, \((0,0,0)\) is an equilibrium with eigenvalues \(\{0,\pm i\}.\) The parameter values \(a=b=c=d=0\)
lead to the well-known Kuramoto--Sivashinsky equation, where it is a divergent free system.
Chang \cite{ChangKS}'s remark implies that the quadratic-truncated classical normal form of this equation falls in \(\LST,\)
\ie a divergent free vector field with a first integral. Here, we also notice that classical normal forms of volume-preserving
Hopf-zero singular systems are not generally a volume-preserving system. Indeed,
a higher degree-truncated normal form of Kuramoto--Sivashinsky system is not divergent free and also does not have a first integral.
Thus, we consider parameter values
\be\label{constant5.9}
b:=3a, \ c:={\frac {7}{975}}\left(31+859a-906a^2\right), \ \hbox{ and } \ d:=-\frac{11}{7}c+10a-10a^2.
\ee Hence, the cubic-truncated first level normal form \(v^{(1)}\) is given by
\ba\nonumber
v^{(1)}&:=& \Theta^0_0+(a-2)F^{-1}_0+ \frac{176a^2+161a-1}{1300} F^0_1+ (1-a)F^1_1
+\frac{1162a^2-1068a-62}{325} F^2_2
\\&&\label{1LevelKS}
-\left(\frac{1013a^2}{3900}-\frac{5737a}{62400}-\frac {14333}{62400}\right)\Theta^0_1
+\frac{1}{2}a\Theta^1_1+\frac{1901a^2-4689a+2399}{1300} \Theta^2_2.
\ea
The infinite level normal form \(v^{(\infty)}\) is formulated by
\be\label{5.12} v^{(\infty)}:=
\left\{
  \begin{array}{ll}
    v^{(\infty)}_+, & \hbox{for } a\in (1,2), \\
    v^{(\infty)}_-, & \hbox{for } a<1 \hbox{ and } a>2,
  \end{array}
\right.
\ee where
\bas\nonumber
v^{\infty}_{\pm}&:=&\Theta^0_0+\frac{1}{2}{F^{-1}_0}\pm F^1_1
\pm \frac a{2(1-a)}\Theta^1_1\pm \frac{\sqrt{2}(15384a^{3}-
67767a^2+94215a-38384)}{-20800(a-1)\sqrt{|(a-2)(a-1)|}} \Theta^2_2
\\ &&
-\frac{\sqrt{2}(a-2)(581a^2-534a-31)}
{325(a-1)\sqrt{|(a-2)(a-1)|}}F^2_2
+\frac{(176a^2+161a-1)(581a^2-534a-31)}{1690000(a-1)^2}F^3_3.
\eas Equation \eqref{5.12} is trivially derived from Proposition \ref{Propcoef}.
\end{exm}

\section{ Optimal truncation } \label{IoosSection}

The convergence analysis of normal forms is an important and difficult problem.
Indeed, the divergent series may appear in both classical normal form computation step and the hypernormalization steps. Generally, the convergence analysis of normal form and transformations have been rarely performed for the hypernormalization steps. Str\'{o}\.{z}yna and \.{Z}oladek \cite{Zoladek02} proved that the normal form series associated with Bogdanov-Takens pre-normal (classical) form are convergent while they \cite{StZolDiv} presented an example of multi-dimensional nilpotent singularity such that the associated normalized (classical) series are divergent.
For the case of divergent normal forms, normal form computations can still be useful. Furthermore, one usually truncates the normal form system at certain grade for the analysis. In this direction we mention possible \emph{jet determinacy} and \emph{optimal truncation} while we only pursue the second one in this section.

The jet determinacy deals with finding a possible degree, say \(k\), for truncation of normal forms such that its dynamics would be qualitatively the same as the untruncated system. This is referred to as \(k\)-determinacy or \(k\)-jet sufficiency; see \cite[Item 4 on Page vi]{MurdBook}. It implies that a convergence analysis for finitely determined systems is not necessary. This is beyond the scope of this paper, but it provides a possible reason (in addition to the involved difficulties) to explain why convergence analysis has not attracted much attention in the literature.

The second idea, that may work for any possible non-finitely determined and divergent normalized system, is as follows. A practical normal form process is stopped at some step and then, one truncates the normal form series at a certain grade. Next, the whole normalized system is considered as a small perturbation of the truncated one. (Note that this perturbation should not be confused with an essentially different concept and materials in Section \ref{PNF}.) This gives rise to an actual first integral for the truncated (unperturbed) system. Therefore, it is important to truncate the normalized system at a grade in which the remainder is (optimally) small. In this direction for any (small) \(\delta\)-neighborhood of the equilibrium, Iooss and Lombardi \cite{IoosLombardi} derived an optimal grade for truncation. They proved that the remainder is analytic and is exponentially small in \(\delta\). Their method works well for a large class of systems that includes Hopf-zero singularity systems. We apply their results in this section to the generalized Kuramoto--Sivashinsky and modified R\"{o}ssler equations in Examples \ref{RosselerExample} and \ref{KSExample}. We restate their results for a simplified and specific case that we need. (There is no claim of novelty in our representation.) For any \(\delta\)-neighborhood, an optimal degree is provided for truncation of the first level normal form systems such that the remainder is exponentially small in \(\delta\).

For any arbitrary natural number \(K,\) the vector \(\lambda=(0, i, -i)\in \mathbb{C}^3\) is called \(\frac{1}{2}, K\)-homologically without small divisors. This is because for any \(\mathbf{m}\in \mathbb{N}^3\) such that \(2\leq |\mathbf{m}|\leq K,\) we have
\bes
\left|\langle\lambda, \mathbf{m}\rangle-\lambda_j\right|\geq \frac{1}{2}, \ \ \hbox{ when } \ \langle\lambda, \mathbf{m}\rangle\neq \lambda_j \ \hbox{ for } \  j= 1, 2, 3,
\ees and \(\lambda=(\lambda_1, \lambda_2, \lambda_3)\); see \cite[Definition 1.2.]{IoosLombardi}.
Now we may present the following proposition.

\begin{prop}\label{IoosTheorem}
Consider an analytic differential system
\ba\label{Eq6.1}
\left\{
  \begin{array}{ll}
\dot{x}= \sum^\infty_{i+j+k=2} a_{i,j,k} x^{i} y^jz^k, &
\\
\dot{y}= z+\sum^\infty_{i+j+k=2} b_{i,j,k} x^{i} y^jz^k,&
\\
\dot{z}= -y+ \sum^\infty_{i+j+k=2} c_{i,j,k} x^{i} y^jz^k. &
  \end{array}
\right.
\ea Let \(\mathfrak{c} \geq 2\) be such that
\be\label{c}
\left\|\sum_{i+j+k=N} a_{i,j,k} x^{i} y^jz^k\ddx+b_{i,j,k} x^{i} y^jz^k\ddy+c_{i,j,k} x^{i} y^jz^k\ddz\right\|\leq \mathfrak{c} \|(x, y, z)\|.
\ee
Then, there exist polynomial transformations such that they send the system \eqref{Eq6.1} into
\ba\label{CNFH0}
\left\{
  \begin{array}{ll}
 \frac{dx}{dt}=
\sum^p_{n=2} a^i_j
x^{n-2i}\rho^{2i}+\mathcal{R}_{x,p}, &
\\
\frac{d\rho}{dt}= \sum^p_{n=2} b^i_jx^{n-2i-1}\rho^{2i+1}+\mathcal{R}_{\rho,p},&
\\
\frac{d\theta}{dt}= 1+ \sum^p_{n=1} c^i_jx^{n-2i}\rho^{2i}+\mathcal{R}_{\theta,p},&
  \end{array}
\right.
\ea
where \(a^i_j, b^i_j,c^i_j\in \mathbb{R}\) and \(\mathcal{R}_p(X)= (\mathcal{R}_{x,p}, \mathcal{R}_{\rho,p}, \mathcal{R}_{\theta,p})\) is an analytic function of \(X=(x, \rho,\theta)\) such that \(\mathcal{R}_p(X)= O(\| X\|^{p+1})\);
see \cite[Equation 1.2]{GazorMokhtari} and \cite[Theorem 1.1]{IoosLombardi}. Furthermore,
for any \(\delta > 0\) such that the optimal degree
\be\label{Popt} p_{\rm opt}:=\left\lfloor\frac{2}{\delta \left({19}\sqrt{3}\mathfrak{c}+6\sqrt{3}\right)e}\right\rfloor \ee satisfies the condition \(p_{\rm opt}\geq 2,\) the remainder
\(\mathcal{R}_{p_{\rm opt}}\) is exponentially small, \ie
\be\label{6.5}
\sup_{\|X\|\leq \delta}\left \|\mathcal{R}_{p_{\rm opt}}(X)\right\|\leq M \delta^2\exp\left(\frac{-2}{\delta e\sqrt{3}({19}\mathfrak{c}+6)}\right),
\ee where
\bas
M&:=&\frac{5\mathfrak{c}}{18}({19}\mathfrak{c}+6)^2\left(\mathfrak{m}\sqrt{\frac{27}{8e}}+4e^2\right)\sim 43.28074575\mathfrak{c} (19\mathfrak{c}+6)^2,
\\
\mathfrak{m}&:=&\sup_{p\in\mathbb{N}}\frac{e^2p!}{p^{p+\frac{1}{2}}e^{-p}}\sim 20.08553692.
\eas
\end{prop}
\bpr The proof trivially follows from \cite[Theorem 1.4. a(ii)]{IoosLombardi} and \cite[Theorem 1.1]{IoosLombardi}.
\epr

Note that the condition \(p_{\rm opt}\geq 2\) is an equivalent condition for \(\delta\leq \frac{1}{\left({19}\sqrt{3}\mathfrak{c}+6\sqrt{3}\right)e}\).
The number
\be\label{ccal} \mathfrak{c}=\sup_{N\geq 2}\left \{\sum_{i+j+k=N} |a_{i,j,k}|+|b_{i,j,k}|+|c_{i,j,k}| \right\},\ee simply satisfies Equation \eqref{c}. We first apply Equation \eqref{ccal} (associated with Proposition \ref{IoosTheorem}) to Examples \ref{RosselerExample} and \ref{KSExample}. Then, we use Lagrange multipliers and a Matlab program in order to compute the least of such numbers \(\mathfrak{c}\) satisfying Equation \eqref{c}. This leads to a comparison with different values of \(\mathfrak{c}\) on the magnitude of the normal form remainder; see Figures \ref{RSM1}--\ref{RSM2} and \ref{EXPKS1}--\ref{EXPKS2}.

In summary, this procedure provides an optimal-degree truncation for the classical normal forms such that the remainder is exponentially small.

\begin{exm}\label{RosselerExample6} Consider the differential system \eqref{RosselerEqn} from Example \ref{RosselerExample} and define
\bas
\mathfrak{c}_2&:=&
  {\frac {2|a(a+1)(d-1)| +|{{a}^{2}-2d}|+2|{a(d{a}^{2}-2)}|+|{a({a}^{2}-4d+2)}|+2|{{a}^{2}-2d{a}^{2}+2}|}{2\sqrt {2-{a}^{2}} |({a}^{2}-2 )| }}
\\&&+\frac{(|a|+1)|d|}{\sqrt {2-{a}^{2}}} +\frac {4|a(2d+1 )|+(|a|+2){\sqrt {2-{a}^{2}}}+4|{d{a}^{2}+1}|+|{a}^{2}+ 4d|+{{a}^{2}}+3|a|+2}{2|{a}^{2}-2|},
\eas
\bas
\mathfrak{c}_3&:=&
\left(\frac{1}{2}\sqrt{2-{a}^{2}}
+\frac{3}{2}(|a|+1)+\frac {(3{a}^{2}+6|a|+3)\sqrt {2-{a}^{2}}+{|a|}^{3}+3{a}^{2}+3|a|+1}{2|{a}^{2}-2|}
\right)|e+r|
\\&&
+\left({\frac {{a}^{4}+5{|a|}^{3}+9{a}^{2}+7|a|+2}{2\sqrt {2-{a}^{2}}\,|{a}^{2}-2 | }}+\frac{3{|a|}^{3}+12a^2+15|a|+6}{2| {a}^{2}-2 |}+{\frac {3{a}^{2}}{2\sqrt {2-{a}^{2}}}}+\frac{ | a |+2}{2}
\right)|e-r|,
\eas
and
\bes
\mathfrak{c}(a):=\max\left\{\mathfrak{c}_2(a), \mathfrak{c}_3(a)\right\}.
\ees
\begin{figure} \label{figRS}
\subfigure[\label{figRSa} \(\mathfrak{c}(a)\) versus \(a\) for \((-0.859,0.875)\).]{\includegraphics[width=.5\columnwidth,height=.2\columnwidth]{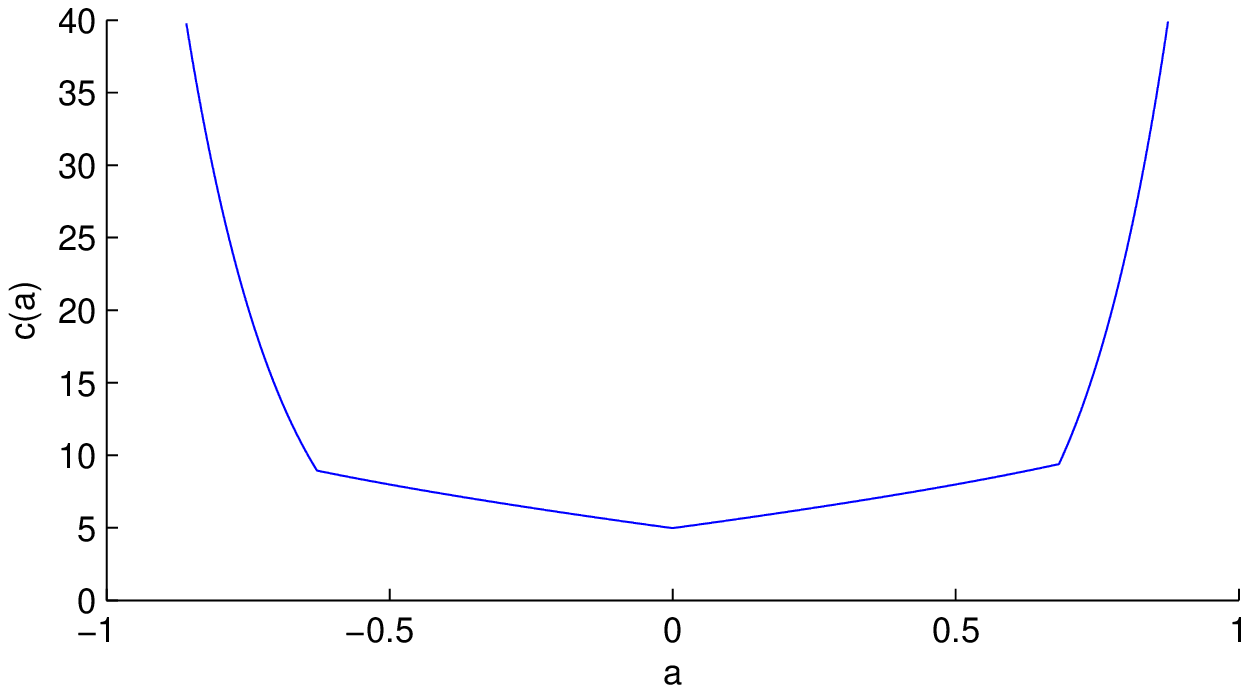}}
\subfigure[\label{figRSb} \(M(a)\) versus \(a\) for \((-0.656, 0.7)\). ]{\includegraphics[width=.5\columnwidth,height=.2\columnwidth]{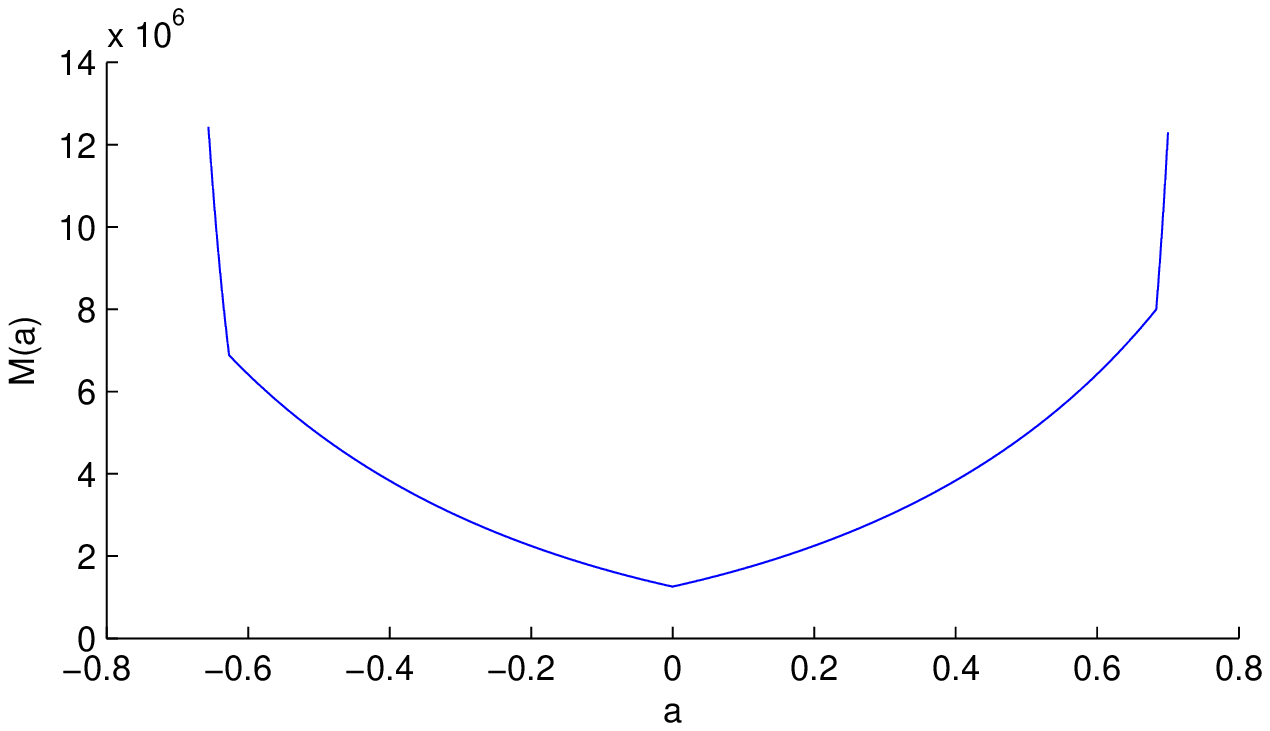}}
\caption{The R\"{o}ssler equation in Example \ref{RosselerExample6}. }\label{caption}
\end{figure}

Recall that \(a\) may take values within the interval \([-\sqrt{2},\sqrt{2}]\). Using a mesh sequence with a step size \(10^{-6}\), the function \(\mathfrak{c}(a)\) is plotted in Figure \ref{figRSa} for \(a\in [-0.859, 0.875]\). Note that the function \(c(a)\) (and hence, \(M(a)\)) is not smooth at \(\pm\frac{1}{\sqrt{2}}\) and the origin. It can be seen that \(\mathfrak{c}(a)\) outside this interval monotonically approaches infinity when \(a\) approaches \(\pm \sqrt{2}.\) The function \(M(a)\) is depicted in Figure \ref{figRSb} using a step size \(10^{-6}\) within the interval \(a\in [-0.656, 0.7]\). Outside this interval, \(M(a)\) monotonically approaches infinity when \(a\) approaches \(\pm \infty.\)

In particular, for \(a=1\) and \(\delta\leq0.6843120638\times 10^{-4},\)
\bas
\mathfrak{c}(1)= 81.52048193,\, M(1)&=&8.530222186\times10^9, \,
p_{opt}(1)=\left\lfloor{ \frac{0.0002731967016}{\delta}}\right\rfloor,
\eas
\be\label{6.6}
\sup_{\| X\|\leq\delta} \|\mathcal{R}_{p_{\rm opt}}(X)\| \leq 8.530222186\times10^9\,{\delta}^{2}{{\rm e}^{-\frac{0.0002731967016}{\delta}
}}.
\ee
Figure \ref{RSM1} depicts the right hand side of Equation \eqref{6.6}.

\begin{figure} \label{RSM}
\subfigure[\label{RSM1} \(\mathfrak{c}\) is obtained via Equation \eqref{ccal}. ]{\includegraphics[width=.5\columnwidth,height=.3\columnwidth]{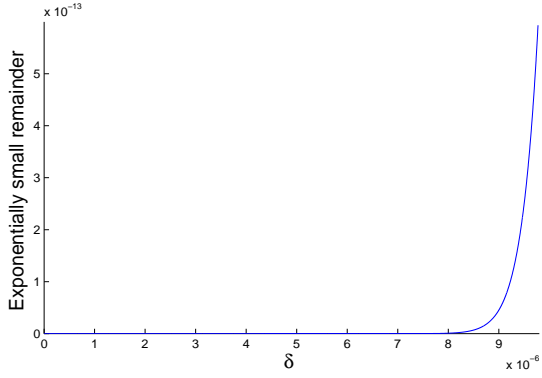}}
\subfigure[\label{RSM2} \(\mathfrak{c}\) is computed by Lagrange multipliers. ]{\includegraphics[width=.5\columnwidth,height=.3\columnwidth]{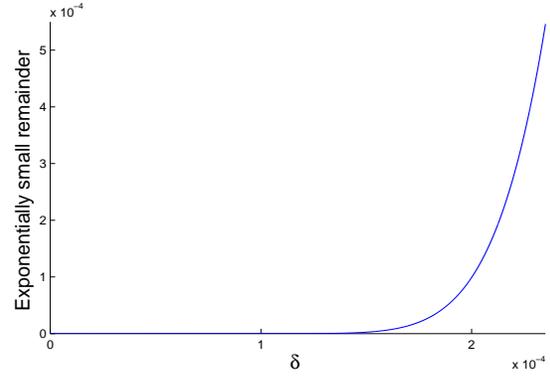}}
\caption{Exponentially small remainder for the R\"{o}ssler equation with \(a=1\).}\label{caption}
\end{figure}

Now using Lagrange multipliers, the least number \(\mathfrak{c}(1)\) satisfying Equation \eqref{c} is given by \(11.72879638\) and the corresponding function associated with the right hand side of Equation \eqref{6.5} is plotted in Figure \ref{RSM2}.
\end{exm}

\begin{exm}\label{KSExample6}
Consider the differential system \eqref{KS} associated with Example \ref{KSExample}. Expanding and simplifying the equations give rise to
\bas
\frac{dx}{dt}&:=&(b-a-2)({z}^{2}+{y}^{2})+2(b+a-2)yz-4x(z+y)-2{x}^{2}
+(c-d)({z}^{3}+{y}^{3})
\\&&
-(c+3d)yz({z}+{y}) -dx(2{z}^{2}+4yz+2{y}^{2}+{x}z+{x}y)
\\
\frac{dy}{dt}&:=&-z+\frac{1}{2}(a-b+2)({z}^{2}+{y}^{2})+(-b -a+2)yz+2xz+2x
y+{x}^{2}+\frac{1}{2}(d-c)({z}^{3}+{y}^{3})
\\&&+\frac{1}{2}(c+{3}d)yz(y{z}+{y})+dx({z}^{2}+2yz+{y}^{2}+\frac{1}{2}{x}z+\frac{1}{2}{x}y)
\\
\frac{dy}{dt}&:=&y+\frac{1}{2}(a-b+2)({z}^{2}+{y}^{2})+(-b -a+2)yz+2xz+2x
y+{x}^{2}+\frac{1}{2}(d-c)({z}^{3}+{y}^{3})
\\&&+\frac{1}{2}(c+{3}d)yz(y{z}+{y})+dx({z}^{2}+2yz+{y}^{2}+\frac{1}{2}{x}z+\frac{1}{2}{x}y)
\eas
We define
\bas
\mathfrak{c}(a):=\max\{ \mathfrak{c}_2(a), \mathfrak{c}_3(a) \},
\eas where
\bas
\mathfrak{c}_2(a)&:=&8|2a-1|+8|a-1|+20,
\\
\mathfrak{c}_3(a)&:=&{\frac {4}{195}} \left|216{a}^{2}+301a-341\right|+{\frac {4}{75}}\left|438{a}^{2}-532a+62\right|+{\frac{4}{325}}\left|2186{a}^{2}-1904a-186\right|,
\eas and \(a\) is an arbitrary parameter. By Proposition \ref{IoosTheorem}, for any
\bes \delta\leq \frac{1}{2e\sqrt{3}({19}\mathfrak{c}+3)} \hbox{ and } p_{\rm opt}:=\left\lfloor\frac{1}{\delta ({19}\sqrt{3}\mathfrak{c}+3\sqrt{3})e}\right\rfloor,\ees
the remainder \(\mathcal{R}_{p_{\rm opt}}\) is exponentially small.

\begin{figure} \label{figKS}
\subfigure[\label{figKSa} \(c(a)\) versus \(a\) for the interval \((-1, 2)\). ]{\includegraphics[width=.5\columnwidth,height=.2\columnwidth]{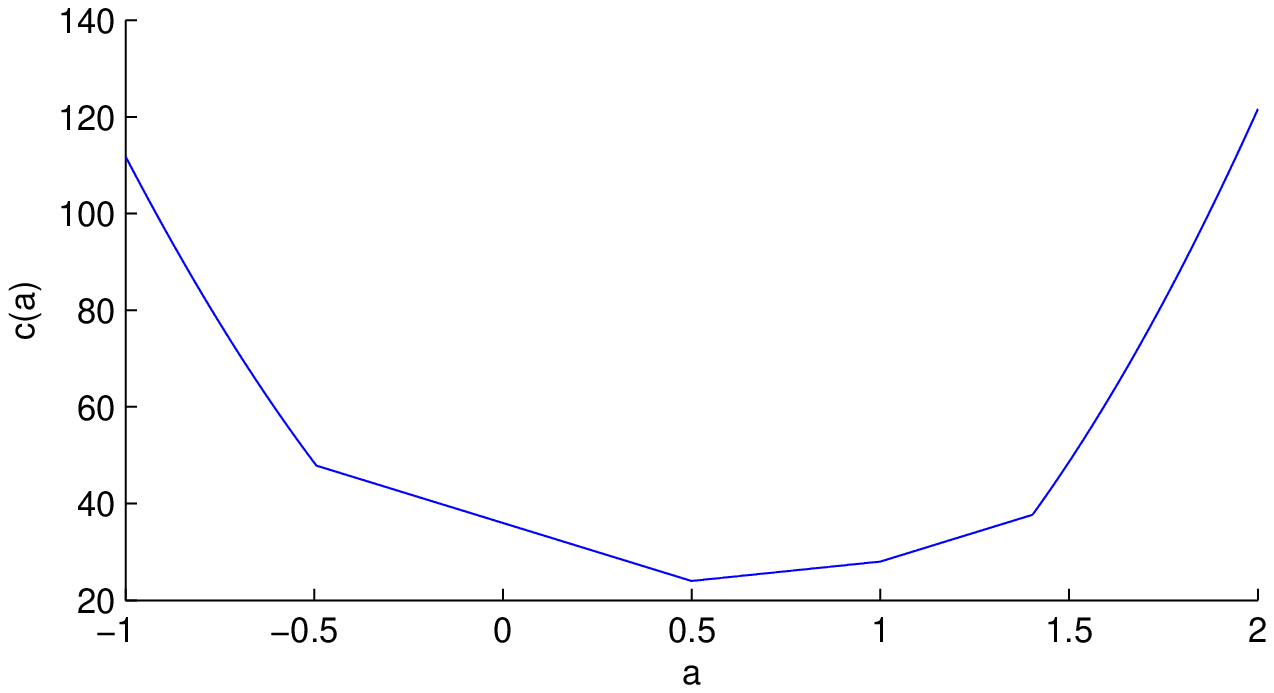}}
\subfigure[\label{figKSb} \(M(a)\) versus \(a\) for \((-0.52,1.515)\). ]{\includegraphics[width=.5\columnwidth,height=.2\columnwidth]{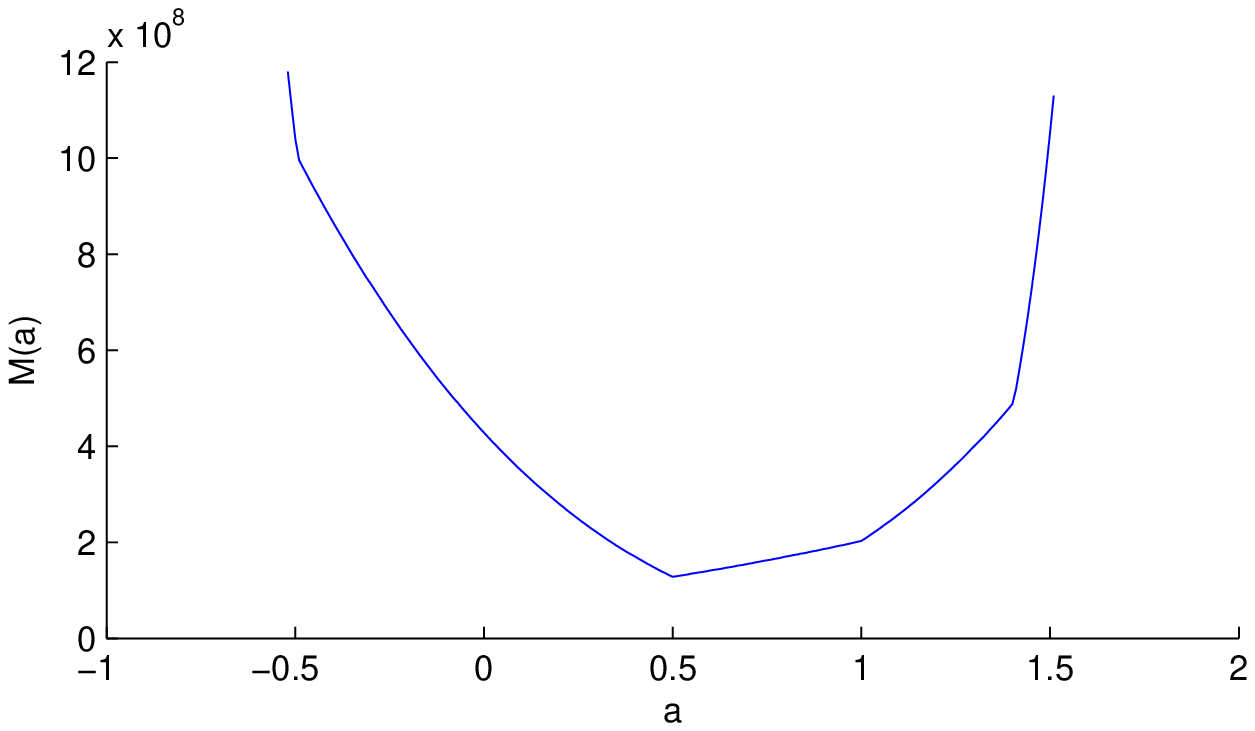}}
\caption{The Kuramoto--Sivashinsky equation in Example \ref{KSExample6}.} \label{caption}
\end{figure}

Figure \ref{figKSa} shows \(\mathfrak{c}(a)\) for \(a\in [-1,2]\) and Figure \ref{figKSb} depicts \(M(a)\) versus \(a\in [-0.52,1.515]\). Outside this interval both functions monotonically approach infinity when \(a\) approaches \(\pm \infty.\) A mesh sequence of step size \(10^{-6}\) is used to plot Figures \ref{figKSa} and \ref{figKSb}. Notice that the function \(\mathfrak{c}(a)\) (and thus, \(M(a)\)) is not smooth at the origin, \(\frac{1}{2}, 1\) (due to \(\mathfrak{c}_2\)), \(-0.494198094\) and \(1.403112056\) (due to the intersection of the curves \(\mathfrak{c}_2\) and \(\mathfrak{c}_3\)).

Now we numerically investigate the exponentially small remainder as a function of \(\delta\). For instance when \(a=1\) for any \(\delta\leq 1.985002751\times 10^{-4},\) we have
\bas
p_{\rm opt}(1)&=&\left\lfloor\frac{0.0003970005502}{\delta}\right\rfloor, \mathfrak{c}(1)= 28, M(1)=3.507658608 \times10^8,
\eas and
\ba\label{Remaind}
\sup_{\| X\|\leq\delta} \|\mathcal{R}_{p_{\rm opt}}(X)\| \leq {\delta}^{2}{{\rm e}^{-\frac{0.0007895735854}{\delta}}}\times3.507658608\times10^8.
\ea
The right hand side of Equation \eqref{Remaind} is drawn in Figure \ref{EXPKS1} versus \(\delta\).

\begin{figure} \label{EXPKS}
\subfigure[\label{EXPKS1} \(\mathfrak{c}\) is computed via Equation \eqref{ccal}.]{\includegraphics[width=.5\columnwidth,height=.3\columnwidth]{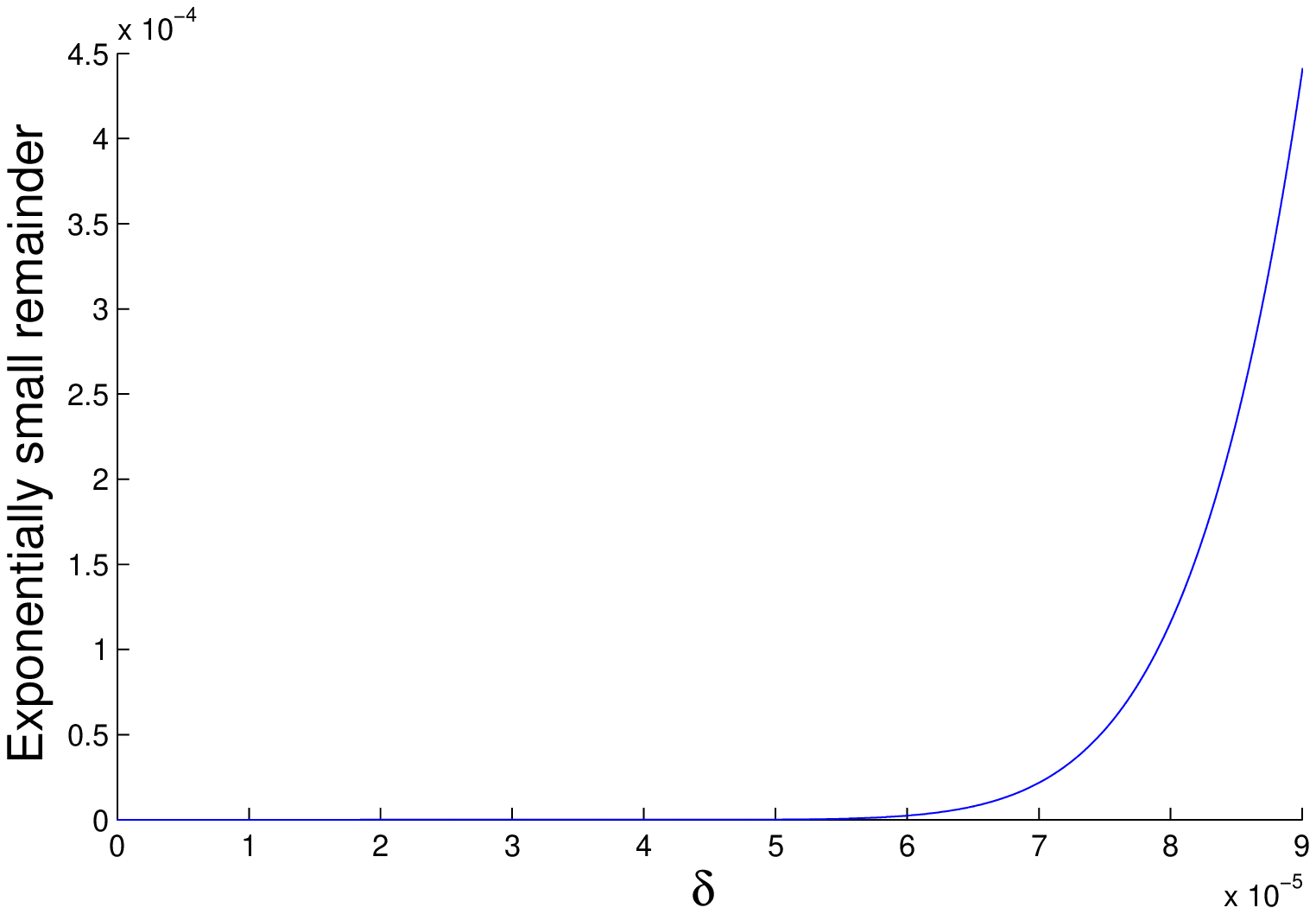}}
\subfigure[\label{EXPKS2} \(\mathfrak{c}\) is obtained by Lagrange multipliers.]{\includegraphics[width=.5\columnwidth,height=.3\columnwidth]{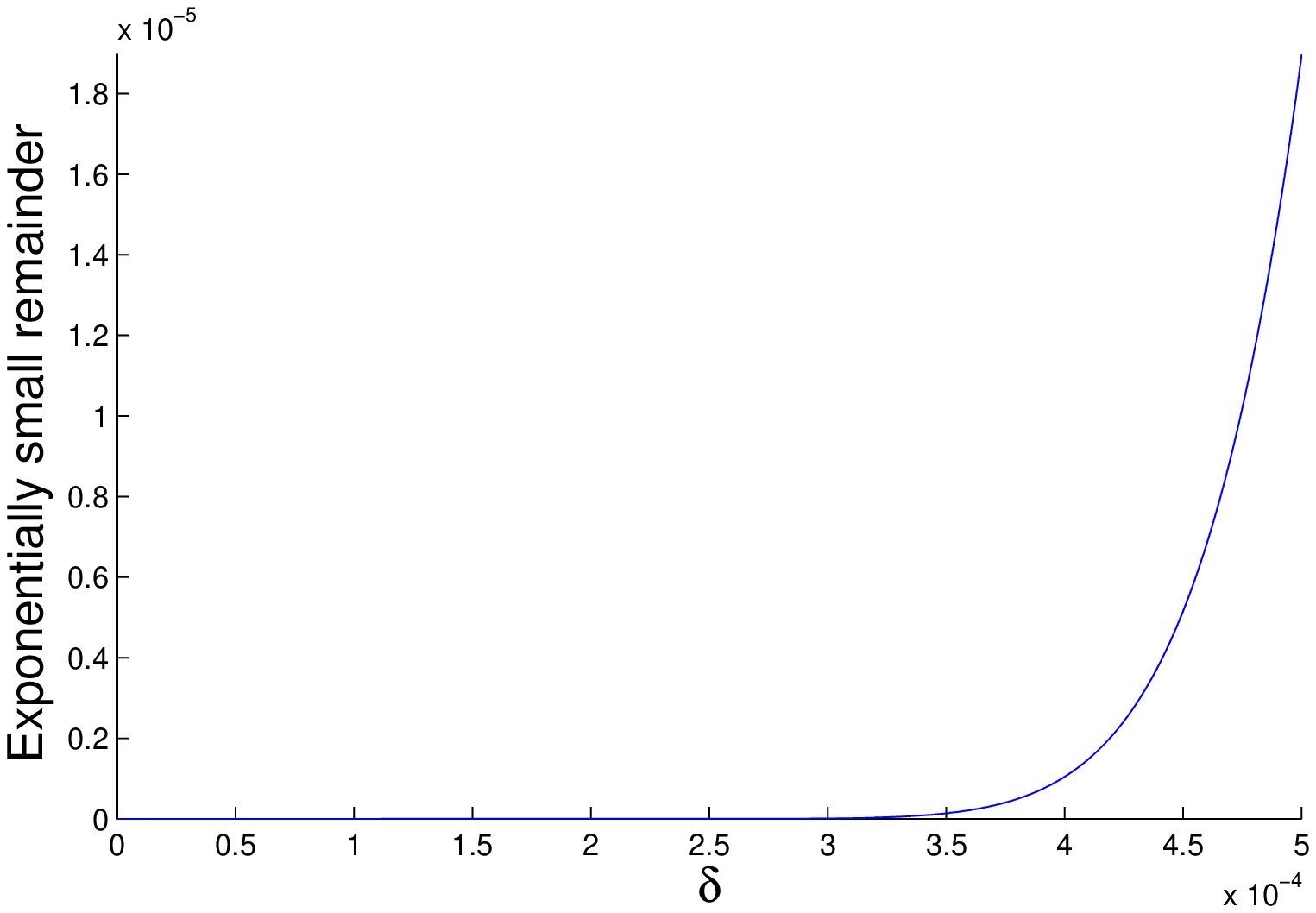}}
\caption{Exponentially small remainder for \(a=1\) in the Kuramoto--Sivashinsky Equation. }\label{caption}
\end{figure}

By using Lagrange multipliers for \(a=1,\) the least number \(\mathfrak{c}\) satisfying Equation \eqref{c} is given by \(4.242640686\).
Accordingly, Figure \ref{EXPKS2} illustrates the changes in the magnitude of upper bound for the remainder associated with \(\mathfrak{c}= 4.242640686\). This compares the nontrivial impact of different values for \(\mathfrak{c}\).

\end{exm}

\section{Radius of convergence }\label{SecConv}

In order to study the convergence of a normalization process two questions arises.
For motivating the first question, we assume that a normal form is computed up to infinite grade as a formal power series. Then, the question is whether or not this power series is convergent.
The second question is with regards to transformations. For this, we recall that the normal form computation is a convergent procedure and the consecutive composition of transformations converges with respect to filtration topology on transformation space; \eg see \cite{GazorYu}. Therefore, one can find a formal power series vector field as a transformation (generator) sending the original system into its normal form up to infinite grade. Hence, the second question is whether the formal transformation series is convergent.

\begin{figure} \label{FigRadios}
\subfigure[\label{RSRadios} R\"{o}ssler Equation ]{\includegraphics[width=.33\columnwidth,height=.33\columnwidth]{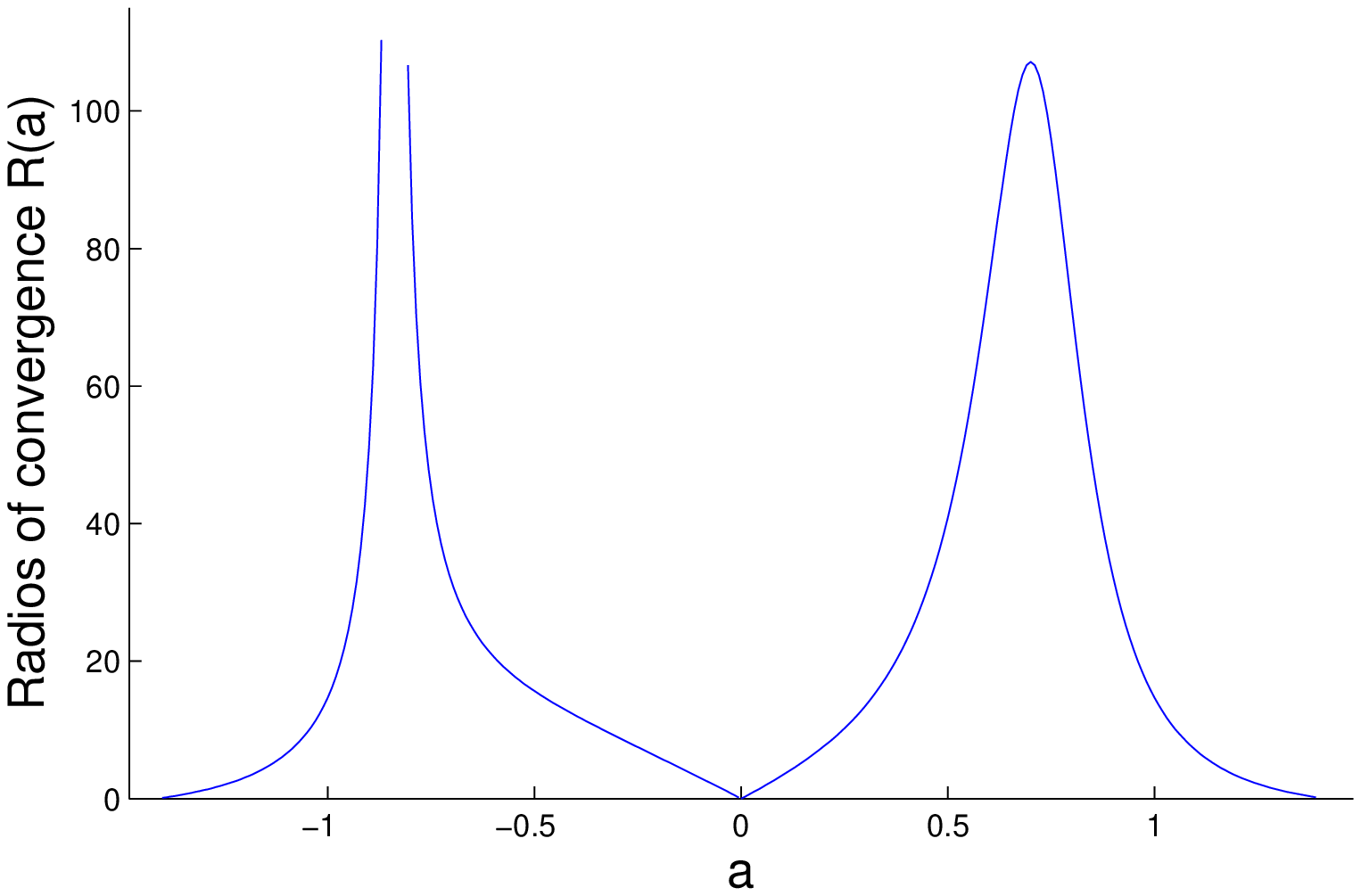}}
\subfigure[\label{KSRadios1} Kuramoto--Sivashinsky]{\includegraphics[width=.32\columnwidth,height=.33\columnwidth]{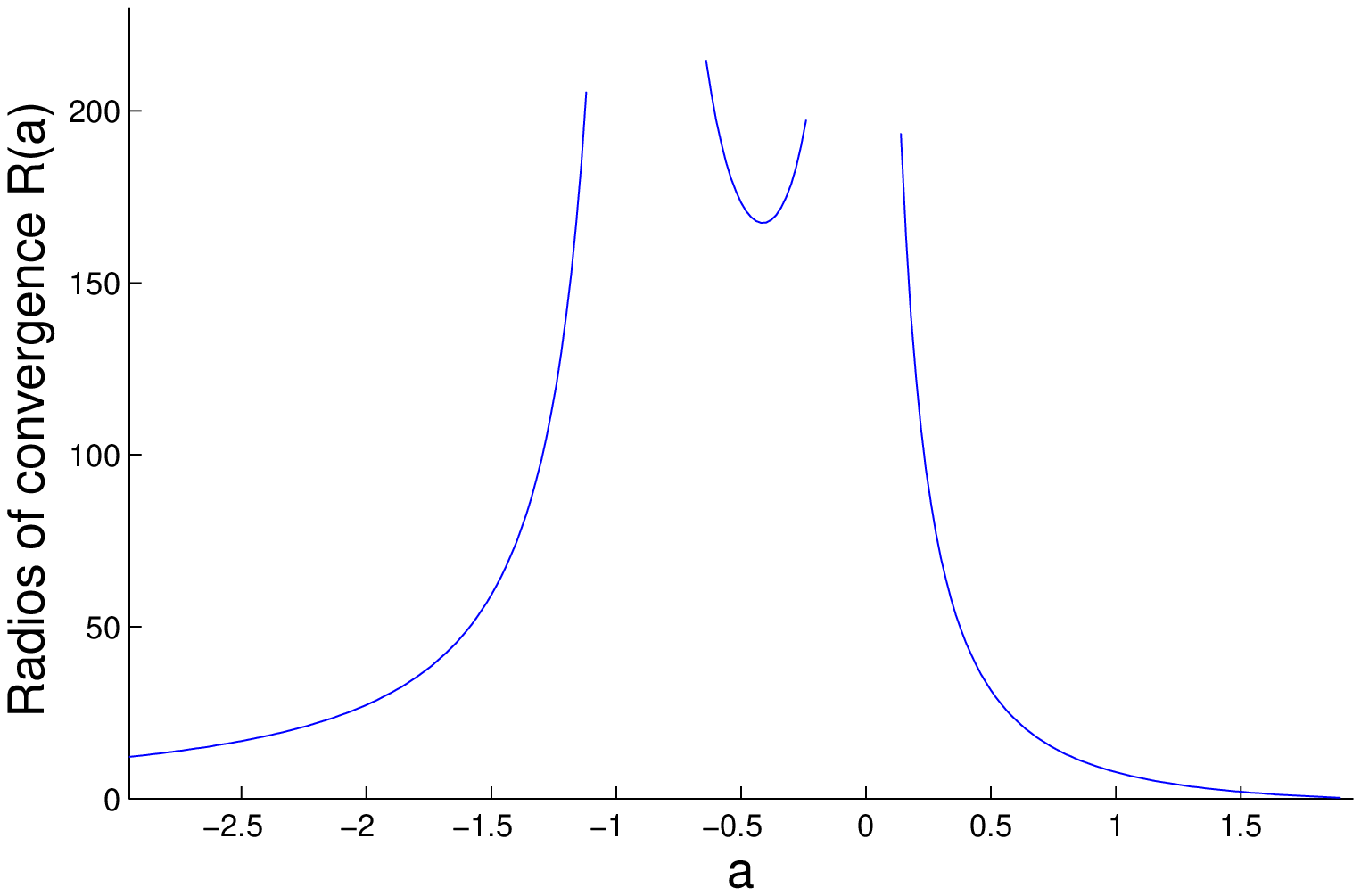}}
\subfigure[\label{KSRadios2} Kuramoto--Sivashinsky]{\includegraphics[width=.33\columnwidth,height=.33\columnwidth]{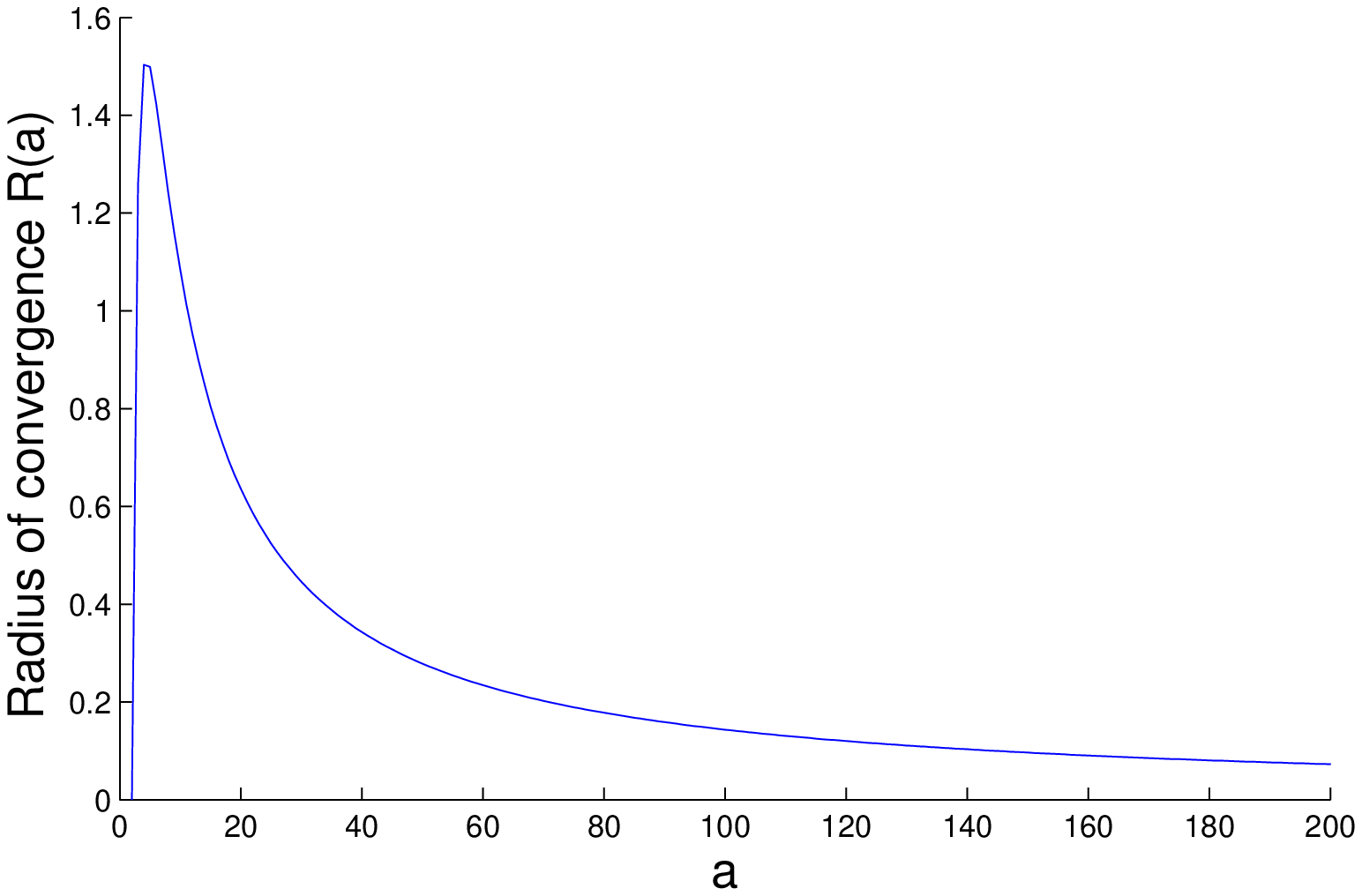}}
\caption{ Numerically suggested radius of convergence. }\label{caption}
\end{figure}

We skip the convergent analysis of transformations in this paper. However, we partially address the first question by computing the normal form of the examples; the modified R\"{o}ssler and generalized Kuramoto--Sivashinsky equations. The second level normal forms of the truncated classical normal form of these examples are computed by Maple up to the grade of one thousand and twenty four. We find the numerically suggested radius of convergence associated with the first integral of the second level normal forms and illustrate them by Figures \ref{RSRadios}, \ref{KSRadios1}, and \ref{KSRadios2}; see Canalis--Durand and Sch\"{a}fke \cite{Gevrey2003,Gevrey2002} where they studied Gevrey type and characters of divergent normal forms.

We numerically analyze a hypernormalization of the truncated first level normal forms of Examples \ref{RosselerExample} and \ref{KSExample} to discuss the convergence of their first integral. The following trivial remark from elementary calculus is the basis of our conclusions in this section. This remark is numbered for keeping the parallel numbering of examples in Sections \ref{Examples}, \ref{IoosSection}, and \ref{SecConv}. This kind of analysis has been rarely performed in the literature; however see \cite{Gevrey2002,Gevrey2003} and the references therein for an advanced numerical convergence analysis of normal forms. The main difficulty rests with an efficient implementation of the results into a computer program. We have implemented the results using Maple XV. The sequence associated with the infinite level normal form does not seem to converge fast. Consequently, our Maple program has not yet been conclusive about radius of convergence for the simplest normal forms.

\begin{rem}
Consider the second level normal form
\begin{eqnarray}\label{Inft}
\dot{x}&=&\rho^2 \pm  x^{p+1}+\sum^\infty_{k=1} \alpha_{k+p} x^{k+p+1},  \\\nonumber
\dot{\rho}&=&\mp\frac{(p+1)}{2} x^{p}\rho-\sum^\infty_{k=1}\frac{(k+p+1)\alpha_{k+p}}{2}x^{k+p}\rho,
\end{eqnarray}
for \(p:=1,\) where \(\alpha_{k+p}\neq 0\) and
\be
L:= \overline{\lim}_{k\rightarrow\infty} \left|\frac{\alpha_{k+p+1}}{\alpha_{k+p}} \right|.
\ee Then, the radius of convergence for the first integral
\bas
f(x, \rho)&:=&\rho^2\left( \frac{1}{2} \rho^2\pm x^{p+1}+\sum^\infty_{k=1}\alpha_{k+p} x^{k+p+1}\right).
\eas is given by \(R:= \frac{1}{L}\).
\end{rem}

\begin{exm} \label{RosselerExample7}
Following Examples \ref{RosselerExample}, we execute our Maple program for Equation \eqref{CubicNFRos} and obtain the second level normal form up to grade one thousand and twenty four for several values of \(a\); that is equivalent to standard degree of five hundred and thirteen. We compute the ratio
\(\alpha_{k+p+1}/\alpha_{k+p}\) and observe that this sequence converges very fast. The numerically suggested radius of convergence for different values of \(a\) is plotted in Figure \ref{RSRadios}. Here, we use a mesh sequence with step-size \(0.01\).

Figure \ref{RSRadios} suggests a critical value at \(a= -0.840563908465308\) whose radius of convergence approaches infinity. This is whence the coefficient of \(F^0_1\) vanishes in the first level normal form. Indeed, the first integral associated with it equals to that of the second level (and also infinite level); see Equation \eqref{F01}.
\end{exm}

\begin{exm} \label{KSExample7}

This example discusses the radius of convergence associated with Example \ref{KSExample}.
By executing our Maple program for Equation \eqref{1LevelKS}, we first obtain the second level normal form up to grade one thousand and twenty four. The ratio \(\alpha_{k+p+1}/\alpha_{k+p}\) is computed for different values of \(a\) and a fast convergence is observed. Therefore, the radius of convergence is approximated for a mesh values of \(a\) with a step-size \(0.01\). The numerically suggested radius of convergence is sketched in Figures \ref{KSRadios1} and \ref{KSRadios2}.

Figure \ref{KSRadios1} suggests two critical values for \(a\); distinct from \(a=2\). These are \(a= -{\frac {161}{352}}\pm {\frac {5}{352}}\sqrt {1065}.\) Here, the coefficient of \(F^0_1\) in the classical level normal form approaches zero; see Equation \eqref{1LevelKS}. Thereby, the first integral associated with Equation \eqref{1LevelKS} is the same as its second and infinite level normal form.

\end{exm}

\end{document}